%% file: ACF.tex
\documentclass[10pt]{amsart}
\usepackage{amssymb}
\usepackage{mathrsfs}
\usepackage{epsfig}

\input xy
\xyoption{all}

\DeclareMathOperator{\Ker}{Ker}

\DeclareMathOperator{\rank}{rank}
\DeclareMathOperator{\dg}{deg}
\DeclareMathOperator{\supp}{supp}
\DeclareMathOperator{\Bbig}{Big}

\DeclareMathOperator{\Amp}{Amp}
\DeclareMathOperator{\Nef}{Nef}

\DeclareMathOperator{\NonDeg}{NonDeg}

\newcommand{\PP}{\ensuremath{\mathbb P}}
\newcommand{\QQ}{\ensuremath{\mathbb Q}}
\newcommand{\RR}{\ensuremath{\mathbb R}}

\newcommand{\CC}{\ensuremath{\mathbb C}}
\newcommand{\OO}{\ensuremath{\mathcal O}}

\newcommand{\rtlsum}[1]{\ensuremath{ \sum_{1\leq l_1< \ldots < l_k\leq m}\hh{#1}{F_{l_1}\cap\dots\cap F_{l_k}}{mpD}}}

\newcommand{\estf}[1]{\ensuremath{\sum_{k=1}^{i+1}{\ #1}}}

\newcommand{\shf}{\ensuremath{{\mathcal F}}}
\newcommand{\sha}{\ensuremath{{\mathcal A}}}
\newcommand{\shb}{\ensuremath{{\mathcal B}}}
\newcommand{\shd}{\ensuremath{{\mathcal D}}}
\newcommand{\shx}{\ensuremath{{\mathcal X}}}

\newcommand{\shl}{\ensuremath{{\mathcal L}}}
\newcommand{\shm}{\ensuremath{{\mathcal M}}}

\newcommand{\hh}[3]{\ensuremath{h^{#1}\left(#2,#3\right)}}
\newcommand{\HH}[3]{\ensuremath{H^{#1}\left(#2,#3\right)}}
\newcommand{\ha}[3]{\ensuremath{\widehat{h}^{#1}\left(#2,#3\right)}}
\newcommand{\sh}[2]{\ensuremath{{\OO}_{#1}\left(#2\right)}}
\newcommand{\euler}[2]{\ensuremath{\chi\left(#1,#2\right)}}
\newcommand{\asy}[1]{\ensuremath{{\mathcal J}\left( \left\| #1 \right\| \right) }}
\newcommand{\ses}[3]{\ensuremath{0\rightarrow #1 \rightarrow #2 \rightarrow #3 \rightarrow 0}}

\newcommand{\vl}[1]{\ensuremath{{\rm vol}\left(#1\right)}}
\newcommand{\vol}[2]{\ensuremath{{\rm vol}_{#1}\left( #2 \right) } }

\newcommand{\zj}[1]{\ensuremath{ \left( #1 \right) }}
\newcommand{\st}[1]{\ensuremath{ \left\{ #1 \right\} }}
\newcommand{\set}[1]{\ensuremath{ \left\{\, #1\, \right\} }}
\newcommand{\norm}[1]{\ensuremath{ \left\| #1 \right\| }}
\newcommand{\deq}{\ensuremath{ \stackrel{\textrm{def}}{=}}}

\theoremstyle{plain}
\newtheorem{theorem}{Theorem}[section]
\newtheorem{lem}[theorem]{Lemma}
\newtheorem*{lemma*}{Lemma}
\newtheorem{lemma}[theorem]{Lemma}
\newtheorem{thm}[theorem]{Theorem}
\newtheorem*{theorem*}{Theorem}
\newtheorem{prop}[theorem]{Proposition}
\newtheorem{proposition}[theorem]{Proposition}
\newtheorem*{proposition*}{Proposition}
\newtheorem{cor}[theorem]{Corollary}
\newtheorem{corollary}[theorem]{Corollary}
\newtheorem*{corollary*}{Corollary}

\theoremstyle{definition}

\newtheorem{remark}[theorem]{Remark}
\newtheorem*{remark*}{Remark}

\newtheorem*{definition*}{Definition}
\newtheorem{definition}[theorem]{Definition}
\newtheorem{ex}[theorem]{Example}
\newtheorem*{example*}{Example}
\newtheorem{example}[theorem]{Example}

\newenvironment{items}
   {\list{\labelitemi}{
      \parsep=0cm \itemsep=0cm \topsep=0cm \partopsep=0.5\baselineskip
      \def\makelabel##1{\hss\llap{\rm##1}}}}
   {\endlist}

\newcommand\newop[2]{\def#1{\mathop{\rm #2}\nolimits}}
\newop\voll{vol}
\newop\NS{NS}
\newop\Neg{Neg}
\newop\Null{Null}
\newop\Pic{Pic}
\newop\Bstab{B_{+}}
\newop\Bst{B_{stab}}
\newop\Bres{B_{restr}}
\newop\B{B}
\newop\Bs{Bs}
\newop\End{End}
\newop\Amp{Amp}
\newop\Face{Face}
\newop\BigCone{Big}
\newop\index{ind}

\newcommand\mif{\ensuremath{\textrm{ if }}}

\newcommand\bbQ{\mathbb Q}

\newcommand\bbR{\mathbb R}

\begin{document}

\title{Asymptotic Cohomological Functions on  Projective Varieties}
\author{Alex K\"uronya}

\maketitle

\bigskip

\section{Introduction}

\input Intro

\section{Basic properties}

\input Asymptotic

\section{Examples}

\input Examples

\section{Cohomological estimates}

\input Estimates

\section{Continuity of Asymptotic Cohomological Functions}

\input Continuity

\section{Auxiliary results}

\input Technical

\medskip

{\small\sc
\noindent
Universit\"at Duisburg-Essen \\
Campus Essen, Fachbereich 6 Mathematik \\
D-45117 Essen, Germany \\
{\em email address:\ }{\tt alex.kueronya@uni-duisburg-essen.de} \\

}

\end{document}

%% file: Intro.tex
Our purpose here is to consider certain cohomological invariants 
associated to complete linear systems on projective varieties. These 
invariants --- called asymptotic cohomological functions --- are 
higher degree analogues of the volume of a divisor.  We establish
the continuity of asymptotic cohomological functions on the real
N\'eron--Severi space and describe  several interesting connections 
which link  them to  classical phenomena, for example
Zariski decompositions of divisors, or Mumford's index theorem for 
the cohomology of line bundles on abelian varieties. 

Our concepts have their origins in the Riemann--Roch problem. 
The classical version asking how $\hh{0}{X}{{\OO}_X(mD)}$ changes as a 
function of $m$  (where $X$ is an irreducible complex  projective variety, and 
$D$ is a Cartier divisor on $X$), has only been answered in dimensions
up to two,  by Riemann and Roch for curves, and by Zariski \cite{Zariski},
and Cutkosky and Srinivas \cite{CutSri} for surfaces. The lack of 
a satisfactory answer in higher dimensions makes it important to look at 
the question from an asymptotic point of view. For ample divisors, the 
by now classical asymptotic Riemann--Roch theorem of Kleiman \cite{Kleiman} 
and  Snapper \cite{Snapper} 
settles the issue. For arbitrary divisors, however, the question
has only surfaced recently in the form of the volume of a divisor, ie. 
 the asymptotic rate of growth of the number of global sections of its multiples. 

The notion  of the volume first arose implicitly in Cutkosky's 
work \cite{cut}, where he used asymptotic computations  to establish 
the non-existence of rational Zariski decompositions on a certain threefold. 
It   was then studied  subsequently  by Demailly, Ein, Fujita, Lazarsfeld,  
and others, while  pioneering efforts regarding other asymptotic 
invariants of linear systems were made by Nakayama \cite{Nakayama}  and Tsuji.
In this
process the  properties of the volume were more fully explored, and 
instead of thinking of  the volume as an invariant linked to a single divisor, 
one  started to consider it as a function on the N\'eron--Severi space, thus as an 
intrinsic  invariant of the underlying variety $X$.

More precisely, the volume of a Cartier divisor $D$ on an irreducible 
projective variety $X$ of dimension $n$ is defined to be 
\[
\vol{X}{D}=\limsup_m{ \frac{ \hh{0}{X}{{\OO}_X(mD)}}{m^n/n!}}\ .
\]
For ample divisors, $\vol{X}{D}= (D^n)$, and 
considered   as a function of the divisor $D$, it descends to a 
degree $n$ homogeneous 
continuous functions  on the real N\'eron--Severi space.
The point of view  that the volume should
be defined on numerical equivalence classes originates in \cite{pag}, 
although it was also realized independently in \cite{Boucksom}.

The volume function is log-concave, which --- according to the 
influential paper \cite{Okounkov} of Okounkov --- is an indication 
that it is a 
'good' notion of multiplicity. In  a different direction, Demailly, 
Ein and Lazarsfeld in \cite{DEL} show that the volume of a divisor 
is the normalized limsup of the moving self-intersection numbers of its 
multiples. The analogous notion  on compact K\"ahler manifolds 
has  been studied by Boucksom \cite{Boucksom}.

Asymptotic cohomological functions of divisors are  direct 
generalizations of the volume function. Let 
$X$ be an irreducible projective variety of dimension $n$, 
$D$ an integral Cartier divisor on $X$. Then for every $0\leq i\leq n$, 
the $i$th {\em asymptotic cohomological
function} associated to $X$ is defined to be
\[
\ha{i}{X}{D} \deq \limsup_{m}{ \frac{\hh{i}{X}{{\OO}_X(mD)}}{m^n/n!}}\ .
\]
Note that by  definition, $\ha{0}{X}{D}=\vol{X}{D}$. There  is a 
pronounced  difference between the case of the volume
function and  higher asymptotic cohomological functions:
for a non-big divisor, the volume is zero, while this is not so
in general for the higher asymptotic cohomological functions. Therefore,  
asymptotic cohomological functions of higher  degree  potentially carry 
information about non-effective divisors as well. 

Probably the most striking property of the volume of divisors  was that 
it defines a continuous function on $N^1(X)_{\RR}$. Our main focus here is 
to prove the corresponding statement for asymptotic cohomological functions. 
More concretely, in Theorem \ref{main}, we establish the following 
(for its $\widehat{h}^0$ predecessor see \cite[Theorem 2.2.44.]{pag}).

\begin{theorem*}[Continuity of asymptotic cohomological functions]
Let $X$ be an irreducible projective variety of dimension $n$ over $\CC$.
Then for all $0\leq i\leq n$, the functions $\widehat{h}^i$ are invariant 
with respect to numerical equivalence of divisors, and 
\[ 
\widehat{h}^i: N^1(X)_{\QQ}\rightarrow {\RR}^{\geq 0}
\] 
defines a continuous  function on $N^1(X)_{\QQ}$. These functions are
homogeneous of degree $n$, 
and satisfy the following Lipschitz-type estimate: there exists a constant $C$ 
such that for all pairs $\xi, \eta\in N^1(X)_{\QQ}$, one has 
\[
| \ha{i}{X}{\xi}-\ha{i}{X}{\eta} |\leq C\cdot \sum_{k=1}^{n}{ \zj{ \max\st{\norm{\xi},\norm{\eta}}}^{n-k}\cdot \norm{\xi-\eta}^k}\ 
\]
for some fixed norm $\|\ \|$.
\end{theorem*}

\begin{corollary*}
With notation as in the Theorem, the asymptotic cohomological 
functions $\widehat{h}^i$ extend uniquely to  continuous functions
\[
\widehat{h}^i: N^1(X)_{\RR}\rightarrow {\RR}^{\geq 0}\ ,
\]
which are homogeneous of degree $n$.
\end{corollary*}

The main ingredients of the proof are boundedness of numerically trivial 
divisors, and cohomological estimates coming from a Mayer--Vietoris-type
exact sequence of sheaves. 

On certain classes of varieties with additional structure, notable 
examples being complex abelian varieties and generalized flag varieties, the 
cohomology of line bundles exhibits an interesting chamber structure, 
in that $\Pic (X)_{\RR}$ is divided into a set of open cones, and for the  
integral points in each such cone there is a single non-vanishing cohomology
group. This behaviour manifests itself in the form of Mumford's index 
theorem for complex abelian varieties, and the Borel--Weil--Bott
theorem on generalized flag varieties. 

These phenomena are suggestively similar to  the behaviour of 
the volume function on smooth projective surfaces. As  described in 
\cite{chambers},  the volume on a smooth projective surface is piecewise 
polynomial with respect to a locally 
finite, locally rational polyhedral chamber decomposition on the cone 
of big divisors. 

Surprisingly enough, there exists some sort of a  generalization of these 
phenomena  to arbitrary irreducible projective varieties, which  rests upon the 
notion of asymptotic cohomological functions. These  
are meaningful on any irreducible projective variety, and in a 
certain sense  they give back the  classical decompositions of 
$\Pic(X)_{\RR}$ on abelian varieties,  and generalized 
flag varieties. 

In many  cases  --- e.g. smooth projective surfaces, toric varieties (see 
\cite{toric})  --- these invariants lead to a chamber 
structure similar to that seen in the case
of abelian varieties or generalized flag varieties, by considering 
the maximal regions in the N\'eron--Severi space where each asymptotic
cohomological function is  given by  a single polynomial.

It is informative to have a look at a concrete example before 
exploring asymptotic cohomological functions in more detail.

\begin{example*}
Consider the projective plane ${\PP}^2$ blown up at a point, which we 
denote by $X$.   
The N\'eron--Severi space of $X$ is generated by the exceptional divisor $E$ 
of the blow-up and  $H$, the pull-back of the class of 
a line in ${\PP}^2$.

The effective  cone is spanned by the rays of $E$ and $H-E$.
A class $\alpha = xH-yE\ \in\ N^1(X)_{\RR}$
is nef in and only if $
\alpha\cdot E\geq 0 \textrm{ and } \alpha\cdot (H-E)\geq 0$, and
the nef cone is generated by the  rays spanned by $H$ and $H-E$. 

\begin{center}
   \epsfig{file=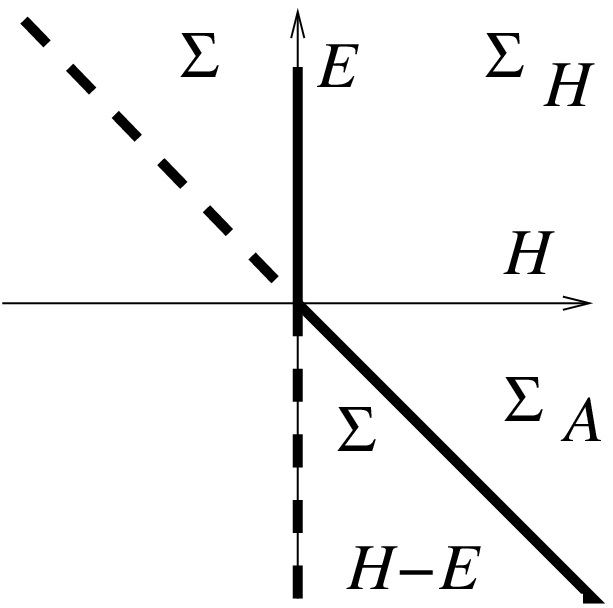,height=1.5in,width=1.5in} \label{fig1}
\end{center}
%\marginnote{Need smaller picture with larger letters}

With the notation of the picture, ${\Sigma}_A$ is the nef cone, the effective
cone is generated by the rays of $E$ and $H-E$, thus it is the union of 
the cones ${\Sigma}_A$ and ${\Sigma}_H$.

We will describe  the asymptotic cohomological functions in some of the 
regions. A  direct computation shows that for $D=xH-yE$, both $x,y>0$
(ie. when $D$ is nef),  
\[
\ha{i}{X}{\alpha}= \begin{cases} {\alpha}^2=x^2-y^2 & \textrm{ if } i=0 \\
                                 0 & \textrm{ if } i=1 \\
                                 0 & \textrm{ if } i=2\ .
                   \end{cases}
\]
The asymptotic cohomological functions are more interesting in the cone 
${\Sigma}_H$ , ie. in the part of $\Bbig(X)$, which consists of non-nef 
divisors. Consider $D=xH+yE$, where $x,y>0$. Then
the asymptotic cohomological functions on ${\Sigma}_L$ are               
\[
\ha{i}{X}{xH+yE}= \begin{cases} x^2 & \textrm{ if } i=0 \\
                                 y^2 & \textrm{ if } i=1 \\
                                 0 & \textrm{ if } i=2\ .
                   \end{cases}
\]
For the part of $N^1(X)_{\RR}$  between the lines of $E$ and $E-H$,
the asymptotic cohomological functions are  given by
\[
\ha{i}{X}{\alpha}= \begin{cases} 0 & \textrm{ if } i=0 \\
                                 y^2-x^2 & \textrm{ if } i=1 \\
                                 0 & \textrm{ if } i=2\ .
                   \end{cases}
\] 
\end{example*}

Results on  ordinary coherent cohomology  have natural  implications for
asymptotic cohomological functions. Among these consequences, one obtains 
 an  asymptotic version of  Serre duality: with notation as above, 
let $D$ be an arbitrary $\QQ$-divisor. Then one has
\[
\ha{i}{X}{D}=\ha{n-i}{X}{-D}\ .
\]

The reason for which we  obtain an  asymptotic duality statement in this 
generality, is  that asymptotic cohomological functions 
are invariant under pullbacks with respect to birational maps. 
More generally, if $f:Y\rightarrow X$ is a finite surjective 
map of   $n$-dimensional irreducible projective 
varieties, $D$ and integral Cartier divisor  on $X$, then 
\[
\ha{i}{Y}{f^*D} = d\cdot \ha{i}{X}{D} \ ,
\]
where  $d$ is the degree of $f$.

Asymptotic cohomological functions yield   a vanishing result, 
which generalizes the asymptotic version of Serre vanishing for 
big divisors.   In particular, a  multiplier ideal argument (Proposition
\ref{basloc}) shows that

\begin{proposition*}
If $X$ is a smooth projective variety over the complex numbers,  
$D$  a big $\QQ$-divisor on $X$ with $d$-dimensional stable base 
locus, then 
\[ \ha{i}{X}{D}=0 \]
for all $i>d$.
\end{proposition*}

Our investigation would not be complete without examples. More than 
 just illustrations of the theory, they were the driving force 
behind many of the developments. We treat abelian varieties, 
homogeneous spaces, and smooth surfaces in detail. In the first two 
cases, we draw on earlier work of Mumford, and Borel--Weil--Bott, 
respectively.

In the case of surfaces, the computation of asymptotic cohomological 
functions  does 
not rely on actual computations of ordinary cohomology groups, rather, 
it makes   use of the  approach in \cite{chambers}, where the authors
 use Zariski  decomposition to determine the volume of divisors.

About the organization of this paper: the basic properties of asymptotic
cohomological functions are given in Section 2. We illustrate these 
functions in a few concrete examples in Section 3. Section 4 hosts the 
technical basis for   this work, it contains  estimates on  
differences of dimensions of cohomology groups of Cartier divisors. 
Based on this, we prove our main result, the continuity of asymptotic 
cohomological functions in Section 5. The last section contains technical 
material, which is used in other parts of the paper, but does  not 
essentially contribute to the understanding of our results. 

%**************************************************************************

\textbf{Acknowledgments.}
I would like to thank  
Tommaso De Fernex, Katalin Friedl, Bill Fulton, Mel Hochster,
Leonardo Mihalcea, Endre Szab\'o and Alex Wolfe for the many 
helpful discussions and comments on previous versions of this 
work.  The observations of Lawrence Ein and Mihnea Popa were particularly 
useful, and influenced this paper considerably. 
I am indebted to Brian Conrad for sharing his technical expertise,
which led to significant improvements.  
Finally, I would like to express my gratitude towards my former 
thesis advisor, Rob Lazarsfeld, for his invaluable guidance
during my graduate years.

%% file: Asymptotic.tex
This section develops  the basic  theory of  asymptotic cohomological
functions  on projective varieties.

First let us  fix notation. In what follows, $X$ will be an irreducible
projective variety over the complex numbers,  unless otherwise mentioned. 
We will use line bundle and divisor notation interchangeably, depending on the 
context. As we will be dealing with $\QQ$- and $\RR$-divisors frequently,
there will be a preference for  divisor language. If it does not
cause confusion, we will use the shorthand notation $\HH{i}{X}{D}$ for  
$\HH{i}{X}{{\OO}_X(D)}$, where $D$ is a Cartier divisor, as it will 
often simplify the appearance of our formulas.

The N\'eron--Severi space $N^1(X)$ is  the group of Cartier 
divisors modulo numerical equivalence. The rational N\'eron--Severi 
space  is denoted by  $N^1(X)_{\QQ}$, 
and similarly for divisor classes with real coefficients. The notations  
$\Nef(X)$ and $\Bbig(X)$ stand for the 
convex cones in $N^1(X)_{\RR}$ 
generated by the classes of nef, and  big divisors, respectively. 

Once past the definition of asymptotic cohomological functions, we will 
present the asymptotic counterparts of some important properties of coherent 
cohomology (e.g. Serre vanishing, Serre duality, K\"unneth formula), and 
describe  the behaviour of asymptotic cohomological functions with respect
to pullbacks. The section ends with a generalization of asymptotic Serre 
vanishing to big divisors. 

\begin{definition}[Asymptotic cohomological functions]
Let $X$ be an irreducible projective variety of dimension $n$, $D$ a 
Cartier divisor  on $X$. The value of the $i$th asymptotic cohomological 
function associated to  $X$ at $D$ is defined to be
\[
\ha {i}{X}{{\OO}_X(D)}\stackrel{\rm\small def}{=}
\limsup_{m}{\frac{\hh{i}{X}{{\OO}_X(mD)}}{m^n/n!}}\ .
\]
\end{definition}

The case $i=0$ is  the volume of the divisor $D$,  and has been studied 
in detail in \cite{chambers},\cite{DEL}, \cite[Section 2.2.C]{pag},
Note that it is established in \cite[Section 11.4.A]{pag},  
that 
\[
\ha{0}{X}{{\OO}_X(D)}=\lim_{m}{\frac{\hh{0}{X}{{\OO}_X(mD)}}{m^n/n!}}\ .
\]
It is not known if the same holds for $i\geq 1$, except in special cases.

\begin{remark}
\label{growth}
Let $X$ be as above,
${\mathcal F}$ a coherent sheaf,
$D$ a Cartier divisor  on $X$. Then there exists a constant $C$ depending on
$X,D$ and ${\mathcal F}$ only, such that 
$\hh{i}{X}{\shf\otimes {\OO}_X(mD)}\leq Cm^n$.
If in addition $D$ is  nef then 
$ \hh{i}{X}{\shf\otimes {\OO}_X(mD)}\leq Cm^{n-1}$ 
for all $i\geq 1$.

Consequently, on the one hand, the value of asymptotic cohomological 
functions is always finite. On the other hand, if $D$ is nef, then  
\[
\ha{i}{X}{{\OO}_X(D)} = 0 
\]
for all $i\geq 1$. This latter statement can be considered as an
asymptotic version of Serre vanishing. In this case, the asymptotic 
Riemann--Roch theorem implies that for a nef divisor $D$, 
\[
\ha{0}{X}{{\OO}_X(D)}=(D^n)\ .
\]
(For proofs of the cited statements see \cite[Section 1.2.B.]{pag}, or 
\cite{RCAV}).
\end{remark}

\begin{example}[Asymptotic cohomological functions on curves]  
As a first example, consider the case of  curves. Let $C$ be
an irreducible projective curve, $L$ a line bundle  on $C$. By  the
asymptotic Riemann--Roch theorem
\[
\hh{0}{X}{L^{\otimes m}}-\hh{1}{X}{L^{\otimes m}} = m\cdot\dg_C L + O(1)\ .
\]
As $L$ is ample if and only if $\deg_C(L)>0$, we obtain that 
\[
\ha{0}{X}{L} = \begin{cases} \deg_C(L) & \mif \deg_C(L)>0 \\
                             0  & \textrm{ otherwise.}
               \end{cases}
\]
Similarly,  
\[
\ha{1}{X}{L} = \begin{cases} 0  & \mif \deg_C(L)\geq 0 \\
                             \deg_C(L)  & \textrm{ otherwise.}
               \end{cases}
\]
\end{example}

\begin{ex}[Weak asymptotic holomorphic Morse inequalities]
In  \cite{Demailly, Demailly2}, Demailly  established a set of inequalities
for the dimensions of the cohomology groups of the difference 
of two nef divisors on smooth varieties. 
These so-called holomorphic Morse inequalities  
can be quickly extended to arbitrary 
projective varieties by appealing to resolution of singularities.  

With notation as so far, let $D=F-G$ be a difference of
two nef Cartier divisors. Then Demailly's  weak holomorphic Morse 
inequality states that
\[
\hh{i}{X}{{\OO}_X(mD)}\leq m^n\frac{F^{n-i}\cdot G^i}{(n-i)!\,i!} + o(m^n)\ ,
\]
hence
\[
\ha{i}{X}{{\OO}_X(D)}\leq \binom{n}{i} F^{n-i}\cdot G^i\ .
\]
\end{ex}

Our next aim is to see  how asymptotic cohomological functions behave under 
'infinitesimal'  perturbations.  Fix an arbitrary
Cartier divisor $D$ and a coherent sheaf $\shf$ on $X$. 
Instead of considering the sequence 
$\hh{i}{X}{{\OO}_X(mD)}$ for $m\geq 1$, we  consider 
$\hh{i}{X}{{\OO}_X(mD)\otimes \shf}$, 
and ask what its properly normalized upper limit is.

\begin{prop}[Invariance under infinitesimal perturbations]
\label{coherent perturbation}
Let $X$ be an $n$-dimensional irreducible projective variety, $\shf$ a
coherent sheaf, $D$ a divisor on $X$. Then
\[
\limsup_{m}{ \frac{\hh{i}{X}{{\OO}_X(mD)\otimes\shf}}{m^n/n!}} =
\rank\, \shf\cdot\ha{i}{X}{{\OO}_X(D)} \ .
\]
\end{prop}

\begin{proof}
We start with the case when $\shf={\OO}_X(N)$ is  an invertible sheaf 
associated to the Cartier divisor $N$.  Let us write $N$ as the 
difference of two effective divisors $N=E-F$, which
do not contain $D$. % \marginnote{What exactly do I need here? Check Rob} 
This way,  we obtain  the following  short exact sequences
\[
\ses{\sh{X}{mD+(E-F)}}{\sh{X}{mD+E}}{\sh{F}{mD+E}} 
\]
\[ 
\ses{\sh{X}{mD}}{\sh{X}{mD+E}}{\sh{E}{mD+E}}\ .
\]
By looking at the corresponding long exact sequences, we see that 
\[
\left| \hh{i}{X}{mD+(E-F)} - \hh{i}{X}{mD+E}   \right|  \leq 
\max{ \st{\hh{i-1}{F}{mD+E}, \hh{i}{F}{mD+E}}}  
\]
and
\[
 \left| \hh{i}{X}{mD+E} - \hh{i}{X}{mD} \right|  \leq  
\max{ \set{ \hh{i-1}{E}{mD+E},\hh{i}{E}{mD+E}}}\ . 
\]
Observe that by Proposition \ref{growth} the right-hand side terms 
in both inequalities  are at most $C'\cdot m^{n-1}$ for 
some positive constant $C'$ independent of $m$. By the triangle 
inequality   we obtain  
\begin{eqnarray*}
 \left| \hh{i}{X}{mD+N} - \hh{i}{X}{mD} \right| &  \leq  & C\cdot m^{n-1}\ 
\end{eqnarray*}
for some constant $C$. After dividing  both sides by $m^n/n!$ and 
taking limsups, we obtain by Lemma \ref{limsup} that
\[ 
\limsup_m{\frac{\hh{i}{X}{mD+N}}{m^n/n!}} = \ha{i}{X}{D} \ .
\]

Next,  consider a direct sum of invertible sheaves 
$\shf=L_1\oplus\ldots\oplus L_r$. Then 
\begin{eqnarray*}
\hh{i}{X}{{\OO}_X(mD)\otimes \shf} 
& = & \sum_{j=1}^{r}{ \hh{i}{X}{{\OO}_X(mD)\otimes L_j}}\ ,
\end{eqnarray*}
therefore,
\begin{eqnarray*}
 \left| \hh{i}{X}{{\OO}_X(mD)\otimes \shf} \right. & - & \left.  r\cdot\hh{i}{X}{{\OO}_X(mD)}
\right|   \\
& \leq & \sum_{j=1}^{r}{ \left| \hh{i}{X}{{\OO}_X(mD)\otimes L_j} -
\hh{i}{X}{{\OO}_X(mD)} \right|} \ ,
\end{eqnarray*}
and consequently, 
\begin{eqnarray*}
\left| \frac{ \hh{i}{X}{ {\OO}_X(mD)\otimes \shf}}{m^n/n!} \right.  & - & \left.  
r\cdot \frac{\hh{i}{X}{{\OO}_X(mD)}}{m^n/n!} \right| \\
&\leq & \sum_{j=1}^{r}{ \left| \frac{ \hh{i}{X}{{\OO}_X(mD)\otimes L_j}}
{m^n/n!} - \frac{\hh{i}{X}{{\OO}_X(mD)}}{m^n/n!} \right| }\ .
\end{eqnarray*}
By the case of line bundles,  the right-hand side of the previous 
inequality converges to zero, hence by (3) of Lemma \ref{limsup} we 
can conclude that
\[
\limsup_{m}{  \frac{ \hh{i}{X}{ {\OO}_X(mD)\otimes \shf}}{m^n/n!}} =
\limsup_{m}{ r\cdot \frac{\hh{i}{X}{{\OO}_X(mD)}}{m^n/n!}} =
 \rank\, \shf \cdot \ha{i}{X}{{\OO}_X(D)} \ .
\]

For the general case, observe that
for every coherent sheaf $\shf$ of rank $r$ on $X$ there exists a map 
$\varphi$,  line bundles $L_1,\dots , L_r$,  and an exact
sequence of coherent sheaves
\[
0 \rightarrow {\mathcal K} \rightarrow {\bigoplus_{j=1}^{r}{L_j}}
\stackrel{\varphi}{\rightarrow} {\shf} \rightarrow {\mathcal C}\rightarrow 0\ ,
\]
where 
$\dim\supp{\mathcal K}\leq n-1 \textrm{ and } \dim\supp{\mathcal C}\leq n-1$.
Split up this sequence into  short exact sequences:
\[
\ses{\mathcal K}{\bigoplus_{j=1}^{r}{L_j}}{\mathcal G}
\textrm{\ \ \   and \ \ \    }
\ses{\mathcal G}{\shf}{\mathcal C}\ ,
\]
where $\mathcal G = \Ker ( \shf\rightarrow \mathcal C)$. It follows that 
the sequences  
\[
\ses{{\OO}_X(mD)\otimes{\mathcal K}}
{\bigoplus_{j=1}^{r}{{\OO}_X(mD)\otimes L_j}}{{\OO}_X(mD)\otimes{\mathcal G}}
\]
and
\[
\ses{{\OO}_X(mD)\otimes{\mathcal G}}{{\OO}_X(mD)\otimes\shf}
{{\OO}_X(mD)\otimes{\mathcal C}}
\]
are then exact.  As 
$\dim\supp {\mathcal K} \leq n-1 \textrm{ and } 
\dim\supp {\mathcal C} \leq n-1$,
we can apply Lemma \ref{support} to both sequences. Hence 
\[
\limsup_{m}{ \frac{ \hh{i}{X}{{\OO}_X(mD)\otimes {\mathcal G}}}{m^n/n!}} 
 =  \limsup_{m}{ \frac{ \hh{i}{X}
{ \oplus_{j=1}^{r}{{\OO}_X(mD)\otimes L_j}}}{m^n/n!}}\ ,
\]
and 
\[
\limsup_{m}{ \frac{ \hh{i}{X}{{\OO}_X(mD)\otimes \shf}}{m^n/n!}}  =  
\limsup_{m}{ \frac{ \hh{i}{X}{{\OO}_X(mD)\otimes {\mathcal G}}}{m^n/n!}}
\ .
\]
Consequently,
\begin{eqnarray*}
\limsup_{m}{ \frac{ \hh{i}{X}{{\OO}_X(mD)\otimes \shf}}{m^n/n!}} & = & 
 \limsup_{m}{ \frac{ \hh{i}{X}
{ \oplus_{j=1}^{r}{{\OO}_X(mD)\otimes L_j}}}{m^n/n!}} \\
& = & \rank\,\shf \cdot \ha{i}{X}{{\OO}_X(D)}\ .
\end{eqnarray*} 
\end{proof}

\begin{corollary}\label{small support}
With notation as above, if $\dim\supp\shf < \dim X$,
then 
\[
\limsup_{m}{ \frac{\hh{i}{X}{{\OO}_X(mD)\otimes\shf}}{m^n/n!}} = 0\ .
\]
\end{corollary}

\begin{prop}[Homogeneity of asymptotic cohomological functions] 
\label{homogeneity}
Let $X$ be an $n$-dimensional irreducible projective variety, 
$D$ a divisor  on $X$, $a>0$ arbitrary integer. Then
\[ 
\ha{i}{X}{aD}=a^n\cdot \ha{i}{X}{D} 
\]
for all $i\geq 0$.
\end{prop}

\begin{proof} 
Homogeneity  is a consequence of  the following statement: if 
\[ 
{\alpha}^{(i)}_{r}\deq \limsup_k{ \frac{\hh{i}{X}{(ak+r)D}}{(ak+r)^n/n!}} \ ,
\]
then
\[ 
{\alpha}_0^{(i)}=\dots ={\alpha}_{a-1}^{(i)}\ .
\]
Grant this for the moment. Then
\[
\ha{i}{X}{aD} =  \limsup_k{ \frac{ \hh{i}{X}{akD}}{k^n/n!}} 
  =   a^n \limsup_{k}{ \frac{ \hh{i}{X}{akD}}{(ak)^n/n!}} 
 =  a^n\ha{i}{X}{D}\ ,
\]
as we wanted to show.

We are left with proving that
${\alpha}_0^{(i)}=\dots ={\alpha}_{a-1}^{(i)}$.
But this follows immediately from  Proposition \ref{coherent perturbation} 
applied  with $N= 2D, \dots , (a-1)D$.
\end{proof}

\begin{remark}
Homogeneity of asymptotic cohomological functions allows us to extend them
to $\QQ$-divisors. For an arbitrary $\QQ$-divisor $D$, set 
\[
\ha{i}{X}{D}\deq \frac{1}{a^n} \ha{i}{X}{aD}\ ,
\]
where $a$ is a  positive integer with  $aD$  integral. 
It follows from Proposition \ref{homogeneity}  that the right-hand side  is 
independent of the choice of $a$.
\end{remark}

Our next goal is  to describe  how asymptotic cohomological functions behave
with respect to pull-backs.

\begin{prop}[Asymptotic cohomological functions  of pullbacks]
\label{cohomology of pullbacks}
Let $f:Y\rightarrow X$ be a proper surjective map of  irreducible projective 
varieties  with $\dim{X}=n$, $D$ a divisor  on $X$. 
\begin{enumerate}
\item If $f$ is  generically finite with $\deg f =d$,  then
\[
\ha{i}{Y}{f^*{\OO}_X(D)} = d\cdot \ha{i}{X}{{\OO}_X(D)}\ .
\]
\item If  $\dim Y > \dim X=n$, then
\[
\ha{i}{Y}{f^*{\OO}_X(D)} = 0 
\]
for all $i$'s. 
\end{enumerate}
\end{prop}

\begin{proof}
Both cases follow from an analysis of  the Leray spectral sequence 
\[
E_2^{pq}(m) = \HH{p}{X}{R^qf_*(f^*{\OO}_X(mD))}\Rightarrow
\HH{p+q}{Y}{f^*{\OO}_X(mD)}\ .
\]
For (1), we show that 
\[
\hh{i}{Y}{f^*{\OO}_X(mD)} = \hh{i}{X}{ (f_*{\OO}_Y)\otimes {\OO}_X(mD)} +
C\cdot m^{n-1}
\]
for a constant $C$ depending on $X,Y,f$ and $D$. The projection formula 
for  direct images  implies that
\[ 
E_2^{pq}(m)=\HH{p}{X}{R^qf_*(f^*{\OO}_X(mD))} \simeq 
\HH{p}{X}{R^qf_*{\OO}_Y\otimes {\OO}_X(mD)}
\]
for all $p,q\geq 0$. As $f$ is generically finite, the higher direct 
image sheaves $R^qf_*{\OO}_Y$  are supported on proper subschemes of $X$.
Hence  by Corollary \ref{small support}
\[
\hh{p}{X}{R^qf_*{\OO}_Y\otimes {\OO}_X(mD)} \leq C\cdot m^{n-1}
\]
for some constant $C$ (not depending on $m$) for all $q$'s except $q=0$.
Therefore, for every $r\geq 2$ and every diagonal $p+q=i$ all terms 
$E_r^{pq}(m)$ --- except
possibly one --- will grow in terms of $m$  at most as 
$C\cdot m^{n-1}$, and the only possible exception will be of the form
\[ 
\hh{i}{X}{ f_*{\OO}_Y\otimes {\OO}_X(mD)} + O(m^{n-1})\ .
\]
Hence
\[
\hh{i}{Y}{ f^*{\OO}_X(mD)} = \hh{i}{X}{ f_*{\OO}_Y\otimes {\OO}_X(mD)} +
O(m^{n-1})\ 
\]
as we wanted. 

As $f$ is generically finite of degree $d$, 
$\rank f_*{\OO}_Y = d$,
therefore  Proposition \ref{coherent perturbation} implies 
\begin{eqnarray*}
\ha{i}{Y}{f^*{\OO}_X(D)} 
& = &  \limsup_{m}{ \frac{ \hh{i}{X}{f_*{\OO}_Y\otimes {\OO}_X(mD)}}
{m^n/n!}} =  d\cdot \ha{i}{X}{{\OO}_X(D)}\ .
\end{eqnarray*}

For the case (2), consider again the  Leray spectral sequence. 
We can see by induction on $r$ that   
$\dim E_r^{pq,(m)} \leq C_r\cdot m^n$ 
for constants $C_r$ independent of $m$. Therefore, 
\[
\hh{i}{Y}{f^*{\OO}_X(mD)} \leq n\cdot C_n \cdot m^n\ ,
\]
from which the proposition follows  as $\dim Y > \dim X =n$.
\end{proof}

\begin{cor}[Birational invariance of asymptotic cohomological functions] 
\label{birational}
Let $f:Y\rightarrow X$ be a proper surjective birational map of irreducible
projective varieties of dimension $n$, let $D$ be a divisor on $X$. Then
\[
\ha{i}{Y}{f^*{\OO}_X(D)} = \ha{i}{X}{{\OO}_X(D)}
\] 
for all  $i\geq 0$.
\end{cor}

\begin{cor}[Asymptotic Serre duality] \label{Serre}
Let $X$ be an irreducible  projective variety of dimension $n$, $D$ a divisor 
  on $X$. Then for every $0\leq i\leq n$
\[
\ha{i}{X}{D}=\ha{n-i}{X}{-D}\ .
\]
\end{cor}
\begin{proof}
Let   $f:Y\rightarrow X$ a resolution of singularities of $X$. 
Then by Corollary \ref{birational} we have
\[ 
\ha{i}{Y}{f^*{\OO}_X(D)}=\ha{i}{X}{{\OO}_X(D)}\ , 
\]
and
\[
\ha{n-i}{Y}{f^*{\OO}_X(-D)}=\ha{n-i}{X}{f^*{\OO}_X(-D)}\ .
\]
Serre duality on the smooth variety $Y$ gives
\[
\hh{i}{Y}{f^*{\OO}_X(mD)}=\hh{n-i}{Y}{K_Y\otimes f^*{\OO}_X(-mD)}
\]
for every $m\geq 1$. By Proposition \ref{coherent perturbation} 
\[
\limsup_m{ \frac{ \hh{n-i}{Y}{K_Y\otimes f^*{\OO}_X(-mD)}}{m^n/n!}} = 
\ha{n-i}{Y}{f^*{\OO}_X(-D)}\ ,
\]
therefore, 
\[
\ha{i}{Y}{f^*{\OO}_X(D)} = \ha{n-i}{Y}{f^*{\OO}_X(-D)}.
\]
This implies
\begin{eqnarray*} 
\ha{i}{X}{{\OO}_X(D)} & = & \ha{n-i}{X}{{\OO}_X(-D)}\ . 
\end{eqnarray*}
\end{proof}

\begin{remark}
By homogeneity, the previous version of Serre duality remains valid for 
$\QQ$-divisors. Also, once having proven the continuity of asymptotic 
cohomological functions on $N^1(X)_{\RR}$ in Corollary \ref{real version}, 
we obtain 
\[
\ha{i}{X}{\xi} = \ha{n-i}{X}{-\xi}
\]
for every  $\xi\in N^1(X)_{\RR}$ and every $0\leq i\leq n$.
\end{remark}

\begin{corollary} \label{topzero}
With notation as above, let $D$ be a big $\QQ$-divisor on $X$.
Then
\[
\ha{n}{X}{D}=0\ .
\]
\end{corollary}
\begin{proof}
As $D$ is big, $-D$ is not, therefore $ \ha{0}{X}{-D}=0$. 
Then  the corollary follows  by the asymptotic version of Serre duality.
\end{proof}

\begin{remark}[K\"unneth formulas for asymptotic cohomological functions]
\label{Kunneth}
Let $X_1, X_2$ be irreducible projective varieties of dimensions 
$n_1$ and $n_2$, $D_1$,$D_2$ Cartier divisors  on  $X_1$ and $X_2$, 
respectively. 
Then 
\[
\hh{i}{X_1\times X_2}{{\pi}_1^*D_1\otimes {\pi}_2^*D_2} 
\,\leq\,
\binom{n_1+n_2}{n_1}\sum_{i=j+k}{ \ha{j}{X_1}{D_1}\cdot 
\ha{k}{X_2}{D_2}}\ .
\] 
for all $i$'s, where ${\pi}_l$ denotes the projection map to $X_l$ ($l=1,2$).
 Furthermore, we have equality in the case $i=0$, which
follows from the observation  that the $\limsup$
in the definition of $\hat{h}^0$ is  a limit 
(\cite[Section 11.4.A.]{pag}).
\end{remark}

In the remainder of this section we discuss  a  connection between asymptotic 
cohomological functions and stable base loci of divisors. As a motivating
example, consider a smooth projective surface $X$ and a big divisor $D$
on $X$.

In Section 3, we will prove the following fact: if $D$ is a big divisor on 
a smooth surface $X$  with Zariski decomposition $D=P_D+N_D$, then 
\[
\ha{i}{X}{D} = \begin{cases} (P_D^2) & \textrm{ if } i=0 \\
                             -(N_D^2) & \textrm{ if } i=1 \\
                             0 & \textrm{ if } i=2\ .
               \end{cases}
\]
As the intersection matrix of $N_D$ is negative definite, $(N_D^2)=0$ if and 
only if $N_D=0$.  On the other hand, as pointed out in \cite{AIBL}, 
$\supp N\subseteq \textrm{B}(D)$.
Therefore, if  $i> \dim\textrm{B}(D)$, the dimension
of the stable base locus of $D$,  then  $\ha{i}{X}{D}=0$. It turns
out, that this  phenomenon  persists on  varieties of higher 
dimension as well.

\begin{prop}[Stable base loci and vanishing of asymptotic cohomological 
functions] \label{basloc}
Let $D$ be a big divisor on a smooth variety $X$, and assume that the 
dimension of the stable base locus of $D$ is $d$. Then 
\[
\ha{i}{X}{{\OO}_X(D)}=0
\]
for all  $i>d$.
\end{prop}
\noindent
\begin{proof}
%\marginnote{Can we do this without multiplier ideals and for general $X$?} 
Let $Z_m\subseteq X$ denote the subscheme defined by the asymptotic multiplier ideal sheaf $\asy{mD}$. Then $\dim{Z_m}\leq d$ for $m$ large. Consider  the exact sequence
\[
0\rightarrow \asy{mD}\otimes {\OO}_{X}\zj{K_X+mD}\rightarrow {\OO}_X\zj{K_X+mD}\rightarrow {\OO}_{Z_m}\zj{K_X+mD}\rightarrow 0\ .
\]
A form of Nadel vanishing for asymptotic multiplier ideals (\cite{pag}, Section 11.2.B) says that for $i>0$ 
\[
\HH{i}{X}{ \asy{mD}\otimes {\OO}_{X}\zj{K_X+mD}}=0\ .
\]
Therefore, if $i>0$ then
\[
\HH{i}{X}{ {\OO}_X\zj{K_X+mD}}=\HH{i}{X}{ {\OO}_{Z_m}\zj{K_X+mD}}\ .
\]
But $ \HH{i}{X}{ {\OO}_{Z_m}\zj{K_X+mD}}=0 $ for $m\gg 0$, 
as $Z_m$ is a $d$-dimensional scheme in this case. Hence 
\[
\hh{i}{X}{{\OO}_X(K_X+mD)}= 0
\]
for $m\gg 0$. By Proposition \ref{coherent perturbation}, this implies
\[
\ha{i}{X}{{\OO}_X(D)} = 
\limsup_m{ \frac{\hh{i}{X}{{\OO}_X(K_X+mD)}}{m^n/n!}} = 0\ ,
\]
as required.
\end{proof}

%% file: Examples.tex
We aim to provide a  pool of examples where  asymptotic 
cohomological functions are  worked out in detail. The ones presented here 
give evidence of  interesting structures in the N\'eron--Severi space 
arising from the functions $\widehat{h}^i$. Some of our examples ---
abelian varieties and generalized flag varieties --- are classical in the 
sense that cohomology groups of line bundles on them have been described
many years ago.

In the case of our other source of examples, smooth surfaces, 
 the computation of asymptotic cohomological functions does
not rely on information about individual cohomology groups of line bundles,
but geometric data in the form of Zariski decompositions. 
The exposition  here  draws heavily on \cite{chambers}. 

Although we do not discuss it here, 
toric varieties form another class where asymptotic cohomological functions
can be  computed explicitly. In addition, they provide 
interesting combinatorial  information (see \cite{toric}).

\begin{remark} Although the examples we cover here  might suggest that 
the asymptotic cohomological functions are piecewise polynomial, this is not 
the case in general. For a counterexample, we refer the reader to 
\cite{chambers}.
\end{remark}

\subsection{Abelian varieties}
 
Let $X$ be a $g$-dimensional complex  abelian variety, expressed  as a 
quotient of  a  $g$-dimensional complex vector space $V$ by a lattice 
$L\subseteq V$. Line bundles on $X$ are  given   in terms 
of Appel--Humbert data, that is,   pairs $(\alpha,H)$,  where
$H$ is a Hermitian form on $V$ such that its imaginary part $E$ 
 is integral on $L\times L$, and 
\begin{eqnarray*}
\alpha & : & L\rightarrow U(1) \textrm{ is a function for which} \\
&& \alpha (l_1+l_2)=\alpha (l_1)\cdot\alpha (l_2)\cdot (-1)^{E(l_1,l_2)}\ .
\end{eqnarray*}

By the Appel--Humbert theorem 
 any  pair  $(\alpha, H)$ determines 
a unique line bundle ${\mathcal L}(\alpha, H)$ on $X$, and every line 
bundle on $X$ is isomorphic to one of the form  ${\mathcal L}(\alpha, H)$ 
for some $(\alpha, H)$. The Hermitian form $H$ on $V$ is the invariant 
two-form associated to $c_1({\mathcal L}(\alpha, H))$.

Fix a line bundle ${\mathcal L}={\mathcal L}(\alpha,H)$.
It is called  {\em nondegenerate}, if $0$ is not 
an eigenvalue of $H$. For a nondegenerate line bundle 
${\mathcal L}$, the number of negative eigenvalues of $H$ is  referred to as 
the {\em index} of ${\mathcal L}$, which we denote  by $\index{\mathcal L}$. 
The case of a positive definite matrix $H$  corresponds to the 
line bundle ${\mathcal L}(\alpha ,H)$ being ample. 

\begin{theorem}[Mumford's index theorem]\cite{Kempf}
With notation as above, let ${\mathcal L}$ be a nondegenerate line 
bundle on $X$. Then
\[
\hh{i}{X}{{\mathcal L}}=\begin{cases} (-1)^{i}{\chi}({\mathcal L})
= \sqrt{ \det_L{E}} & \mif  i=\index (\mathcal L) \\ 
0 & {\rm\ otherwise.} \end{cases}
\]
\end{theorem}
The first Chern class of  line bundles is additive, therefore
$\index ({\mathcal L}^{\otimes m}) = \index ({\mathcal L})$ for all $m\geq 1$.

Consider $\NonDeg (X)\subseteq \Pic(X)_{\RR}$,
the cone generated by all  nondegenerate line bundles.
Then    $\NonDeg (X)$ is an open cone  in 
$\Pic (X)_{\RR}$, and its complement has Lebesgue measure zero. For every 
$1\leq j \leq g$, define ${\mathcal C}_j$ to be the 
cone spanned by nondegenerate line bundles of index $j$. Then (apart from
the origin) $\NonDeg (X)$ is the disjoint union of ${\mathcal C}_1, \dots 
, {\mathcal C}_g$. On each ${\mathcal C}_j$, we have 
\[
\ha{i}{X}{\xi}=\begin{cases}(-1)^{i} (\xi^n) & {\rm\ if\ } i=j \\ 0 & 
{\rm\ otherwise.} \end{cases}
\]
This way, we obtain 
a finite decomposition of $N^1(X)_{\RR}$ into a set of cones, such that 
on each cone, the asymptotic cohomological functions are homogeneous 
polynomials of degree  $g$.

\begin{example}
We will consider the following example in more detail. Consider  
$X \deq E\times E$, the product of an elliptic curve with itself.  
Fix a point $P\in E$. Then the three classes
\[ 
e_1 = [ \{P\}\times E]\ ,\ e_2 = [E\times \{P\}] \ ,\ \delta = [\Delta]
\]
in $N^1(X)_{\RR}$ are independent ($\Delta\subseteq E\times E$ is the diagonal) and generate $N^1(X)_{\RR}$. The various intersection numbers among them are as follows:
\[
\delta\cdot e_1 = \delta \cdot e_2 = e_1\cdot e_2 = 1
\textrm{ and }
e_1^2=e_2^2 = \delta^2 = 0\ .
\]
Any effective curve on $X$ is nef, $\overline{\textrm{NE}}(X)=\Nef(X)$,
furthermore, a class $\alpha\in N^1(X)_{\RR}$ is nef if and only if
$\alpha^2\geq 0\  \textrm{ and } \alpha\cdot h\geq 0$ 
for some (any) ample class $h$. 

In particular, if $\alpha = x\cdot e_1 + y\cdot e_2 + z\cdot \delta$
then $\alpha $ is nef if and only if 
\begin{eqnarray*}
xy+xz+yz & \geq & 0 \textrm{ and } x+y+z  \geq 0\ .
\end{eqnarray*}
As a reference for these statements see \cite{pag}, Section 1.5.B.
One can see  that the $\Nef(X)$ is a circular cone inside 
$N^1(X)_{\RR}\simeq {\RR}^3$. By continuity, define the index of a 
real divisor class $\alpha\in \NonDeg(X)$ to be $0$ if it is ample, 
$2$ if $-\alpha$ is ample, and $1$ otherwise. Then 
\[
\ha{i}{X}{\alpha}= \begin{cases} (-1)^{\index(\alpha)}(\alpha^2) = 
(-1)^{\index(\alpha)}(xy+xz+yz) & \mif i= \index(\alpha) \\
                                 0 & \textrm{otherwise}.
                   \end{cases}              
\]
\end{example}

%**************************************************************************

\subsection{Smooth surfaces}
For a smooth projective surface $X$ over $\CC$,  we 
exhibit a locally finite  decomposition of 
$\Bbig(X)\subseteq N^1(X)_{\RR}$  into locally polyhedral cones, 
such that on each chamber of the decomposition, the functions 
$\widehat{h}^i$ are given by homogeneous quadratic polynomials
coming from intersection numbers of divisors on $X$. 

The discussion is based on \cite{chambers}, where the
authors work out theory of the 
volume function for surfaces,  however,  the issue of asymptotic 
cohomological functions 
was not raised. Here we treat the general case building on their 
results. Also, we will make use of the continuity of asymptotic 
cohomological functions, which we prove in Section 5.

First  recall some important facts about our main tool, 
Zariski decompositions.

\begin{theorem}[Existence and uniqueness of Zariski decompositions for $\RR$-di\-vi\-sors, \cite{KMM87}, Theorem 7.3.1]
   Let $D$ be a pseudo-effective $\bbR$-divisor on a smooth
   projective surface. Then there exists a unique effective
   $\bbR$-divisor \[ N_D=\sum_{i=1}^m a_iN_i \]  such that
   \begin{items}
      \item[(i)]
         $P_D=D-N_D$ is nef,
      \item[(ii)]
         $N_D$ is either zero or its intersection matrix
         $(N_i\cdot N_j)$ is negative definite,
      \item[(iii)] 
         $P_D\cdot N_i=0$ for $i=1,\dots, m$.
   \end{items}
   Furthermore, $N_D$ is uniquely determined by
   the numerical equivalence class of $D$, and
   if $D$ is a $\bbQ$-divisor, then so are $P_D$ and $N_D$.
   The decomposition 
   \[
   D=P_D+N_D
   \]
   is called the {\em Zariski decomposition} of $D$.
\end{theorem}

The connection between Zariski decompositions and 
asymptotic cohomological functions comes from the following result.

\begin{proposition}[Section 2.3.C., \cite{pag}]
   Let $D$ be a big integral divisor, $D=P_D+N_D$ the Zariski decomposition of $D$. Then
   \begin{items}
      \item[(i)] $\HH{0}{X}{kD}=\HH{0}{X}{kP_D}$ for all $k\geq 1$ such that $kP_D$ is integral, and
      \item[(ii)] $ \vl{D}=\vl{P_D}=\zj{P_D^2}.$
   \end{items}
\end{proposition}

By  homogeneity and continuity of the volume we obtain that
for an arbitrary big $\RR$-divisor $D$ with Zariski decomposition 
$D=P_D+N_D$ we have $  {\rm vol}(D)=\zj{P_D^2}=\zj{D-N_D}^2$.

Let $D$ be an $\RR$-divisor on $X$. In determining the 
asymptotic cohomological functions on $X$, we distinguish 
three  cases,  according to whether $D$ 
is pseudo-effective, $-D$ is pseudo-effective or none. 

\begin{prop}\label{asymptotic for pseff}
With notation as above, if $D$ is pseudo-effective then
\[
\ha{i}{X}{D} = \begin{cases} (P_D^2) & \textrm{ if } i=0 \\
                             -(N_D^2) & \textrm{ if } i=1 \\
                             0 & \textrm{ if } i=2\ .
               \end{cases}
\]
\end{prop}

\begin{proof}   
If $D=P_D+N_D$ is the Zariski decomposition of the 
pseudo-effective   divisor $D$, then $\ha{0}{X}{D}=\zj{P_D^2}$.
Furthermore,   if $D$ is pseudo-effective then by  Corollary \ref{topzero}  
and the continuity of $\widehat{h}^2$, $\ha{2}{X}{D}=0$.
In order to compute $\hat{h}^1$, consider the equality
\[
\hh{1}{X}{mD}=\hh{0}{X}{mD}+\hh{2}{X}{mD}-\euler{X}{mD}\ .
\]
This implies that 
\[
\ha{1}{X}{D}=\limsup_m{ \zj{
\frac{\hh{0}{X}{mD}}{m^2/2}+\frac{\hh{2}{X}{mD}}{m^2/2}
-\frac{\euler{X}{mD}}{m^2/2} } }\ .
\]
All three sequences on the right-hand side are convergent. The $h^0$
sequence by the fact that the volume function is in general a limit. The 
$h^2$ sequence converges by  $\ha{2}{X}{D}=0$. Finally, the convergence 
of the  sequence of Euler  characteristics follows from  the 
Asymptotic Riemann--Roch theorem. 
Therefore the $\limsup$ on the right-hand side is  a limit, and 
$\ha{1}{X}{D}=  - (N_D^2)$.
\end{proof}

\begin{corollary} \label{asymptotic for -pseff and for none}
If $-D$ is pseudo-effective with Zariski decomposition $-D=P_D+N_D$ then
\[
\ha{i}{X}{D} = \begin{cases} 0 & \textrm{ if } i=0 \\
                             -(N_{-D}^2) & \textrm{ if } i=1 \\
                             (P_{-D}^2) & \textrm{ if } i=2\ .
               \end{cases}              
\]
When neither $D$ nor  $-D$ are pseudo-effective,  one has 
\[ 
\ha{i}{X}{D}= \begin{cases} 0 & \textrm { if } i=0 \\
                            -(D^2) & \textrm{ if } i=1 \\
                            0 & \textrm{ if } i=2 \ .
              \end{cases}              
\]
\end{corollary}

As in \cite{chambers}, with a careful examination of the variation of 
Zariski decompositions, one
can  give a geometric description of the volume, and hence all asymptotic 
cohomological functions on the big cone of  a smooth surface. 
The main result is the following. 

\begin{theorem}
With notation as above, there exists a   locally finite decomposition of 
 $\Bbig(X)$ into rational locally polyhedral subcones such that on each of 
those the asymptotic cohomological functions are given by a single 
homogeneous quadratic polynomial.
\end{theorem}

\begin{proof} Follows from our description of asymptotic cohomological 
functions on smooth projective surfaces, and the main theorem of 
\cite{chambers}.
\end{proof}

In some cases, the locally finite chamber structure will turn out to be 
finite polyhedral. 

\begin{prop}
Let $X$ be a del Pezzo surface. Then there exists a finite decomposition of $N^1(X)_{\RR}$ into rational polyhedral cones such that on each of these cones all asymptotic cohomological functions are given by homogeneous quadratic polynomials.
\end{prop}
\begin{proof}
The statement is proved by considering the effective cone, its negative and 
the remaining part separately. Of these three the first two are convex 
rational polyhedral cones, the third one is not convex, nevertheless
 its finitely many boundary components are still rational polyhedral. 
By Proposition \ref{asymptotic for pseff} and 
\cite[Proposition 3.4.]{chambers}, 
the statement of the corollary holds for all asymptotic cohomological 
functions on the effective cone. 
Analogously, Corollary \ref{asymptotic for -pseff and for none}  
implies the same on the negative of the effective cone, 
and  for all classes $\alpha$ where neither 
$\alpha$ nor $-\alpha$ is effective.
\end{proof}

\subsection{Generalized flag varieties}
Let $G$  denote a simply-connected semisimple complex Lie group, 
$B\subseteq G$ a Borel subgroup,  $\Delta,{\Delta}_+$  the set of roots 
and positive roots, respectively. The factor space $X=G/B$ is 
equipped with the structure of an irreducible  projective variety over $\CC$, 
and  there is a  natural isomorphism 
\[
{\Lambda}_{W}\simeq \text{Pic}(G/B)
\]
given by 
$\lambda \mapsto L_{\lambda}=G\times_B {\CC}_{\lambda} $,
where ${\Lambda}_W$ is the associated weight lattice 
(cf. \cite{fulhar}, Section 23.3.). Set 
\[
\rho \deq  \frac{1}{2}\sum_{v\in {\Delta}_+}{v}\ .
\]
The computation of the cohomology groups of the line bundles 
$L_{\lambda}$ is a   celebrated  result of Borel--Weil and Bott.
\begin{thm}[Borel--Weil--Bott, \cite{bott}]
For a given weight $\lambda$ the following (mutually exclusive) situations can happen.
\begin{enumerate}
\item $\lambda+\rho$ is on the boundary of a fundamental chamber of $W$ in ${\Lambda}_W$. Then 
\[ \HH{i}{G/B}{L_{\lambda}}=0 \] for all $0\leq i\leq \dim{G/B}\ $.
\item $\lambda+\rho$ is in the interior of a chamber. Then there is a unique element $w\in W$ such that $w(\lambda+\rho)$ is in the positive chamber. In this case
\[ \HH{i}{G/B}{L_{\lambda}}=\begin{cases}  \HH{0}{G/B} {L_{w(\lambda+\rho)} } & \text{ if $i={\rm ind}(\lambda)$} \\
                                        0 &  \text{otherwise} \ .
                                \end{cases}
\]
\end{enumerate}
\end{thm}
\noindent
Here ${\rm ind}(\lambda)$ is the number of positive roots $v$ such that  $(\lambda+\rho.v)<0$ with respect to the Killing form.
This gives rise to  a decomposition of $N^1(G/B)_{\RR}$ 
into finitely many open polyhedral chambers, on each of which 
every  $\ha{i}{G/B}{L}$ is given by a single homogeneous polynomial. 
Let us describe this chamber structure in more detail. 

To every $v\in {\Delta}_+$ we attach the the half-spaces
\[
H_{v,\rho}^{\pm}=
\left\{ \lambda\in {\Lambda}_W | (\lambda+\rho,v)\gtrless 0 \right\}
\textrm{ and }
H_{v}^{\pm}=\left\{ \alpha\in N^1(G/B)_{\RR} | (\alpha,v)\gtrless 0 \right\}\ .
\]
For every choice of signs, the intersection 
${H_{v_1,\rho}^{\pm}}\cap \dots {H_{v_n,\rho}^{\pm}}$
 --- where $v_1,\dots ,v_n$ are {\em all} the positive roots --- is  either 
empty or an open polyhedral cone, we will also refer to them as 
cohomology chambers.

On such a chamber ${\rm ind}(\lambda)$ is constant, hence by the Borel--Weil--Bott theorem
\[
\hh{i}{G/B}{L_{\lambda}}=\begin{cases}  \hh{0}{G/B} {L_{w(\lambda+\rho)} }  & \text{ if $i=$ind$(\lambda)$} \\
                                0 & \text{otherwise.}
                        \end{cases}
\]
If  $I\subseteq {\Delta}_{+}$ is  an arbitrary subset of positive roots, 
define $C_I$ to be the set of weights $\lambda\in {\Lambda}_W$
for which the sequences $(m\lambda+\rho ,v)$
are eventually positive  (for $m\gg 0$) if and only if $v\in I$. Then
\[
C_I = {\Lambda}_W\cap \bigcap_{v\in I}{ H_{v}^+}  
\cap \bigcap_{v\not\in I}{H_{v}^-}\ .
\]

\begin{corollary}
With notation as above
\[
\ha{i}{G/B}{\alpha}=\begin{cases} (-1)^i ({\alpha}^n) & 
                                \mif  i=\textrm{ind}(\alpha) \\
                                0 & \textrm{ otherwise}.
                        \end{cases}
\]
where $\alpha\in N^1(X)_{\RR}$ and $\index (\alpha)$ is defined to be the number of positive roots $v$ for which $(\alpha, v)>0$ holds.
\end{corollary}

The hyperplanes $H_v$ (where $v$ runs through ${\Delta}_+$) 
determine a finite rational polyhedral decomposition of $N^1(X)_{\RR}$, such 
that on each piece, the asymptotic cohomological functions are given by 
homogeneous polynomials. 
The self-intersection  $({\alpha})^n$ 
can be determined from the Weyl dimension formula 
(see  \cite{fulhar}, Section 24.1).

Let us illustrate  the previous discussion on a concrete example.
\begin{example}
Let $G=SL(3,\CC)$. Then  the upper triangular matrices in $G$ form a 
Borel subgroup $B$. Consider  the root system $A_2$ attached to 
$SL(3,\CC)$, let us denote the three positive roots by $v_1, v_2, v_3 $. 
Assuming that the root vectors have unit length, they will be
$v_1=(1,0), v_2=(\tfrac{1}{2},\tfrac{\sqrt{3}}{2}),  
v_3=(-\tfrac{1}{2},\tfrac{\sqrt{3}}{2})$.
The cohomology chambers  and the asymptotic cohomology chambers (in that
order)  look as follows.

\begin{center}
\epsfig{file=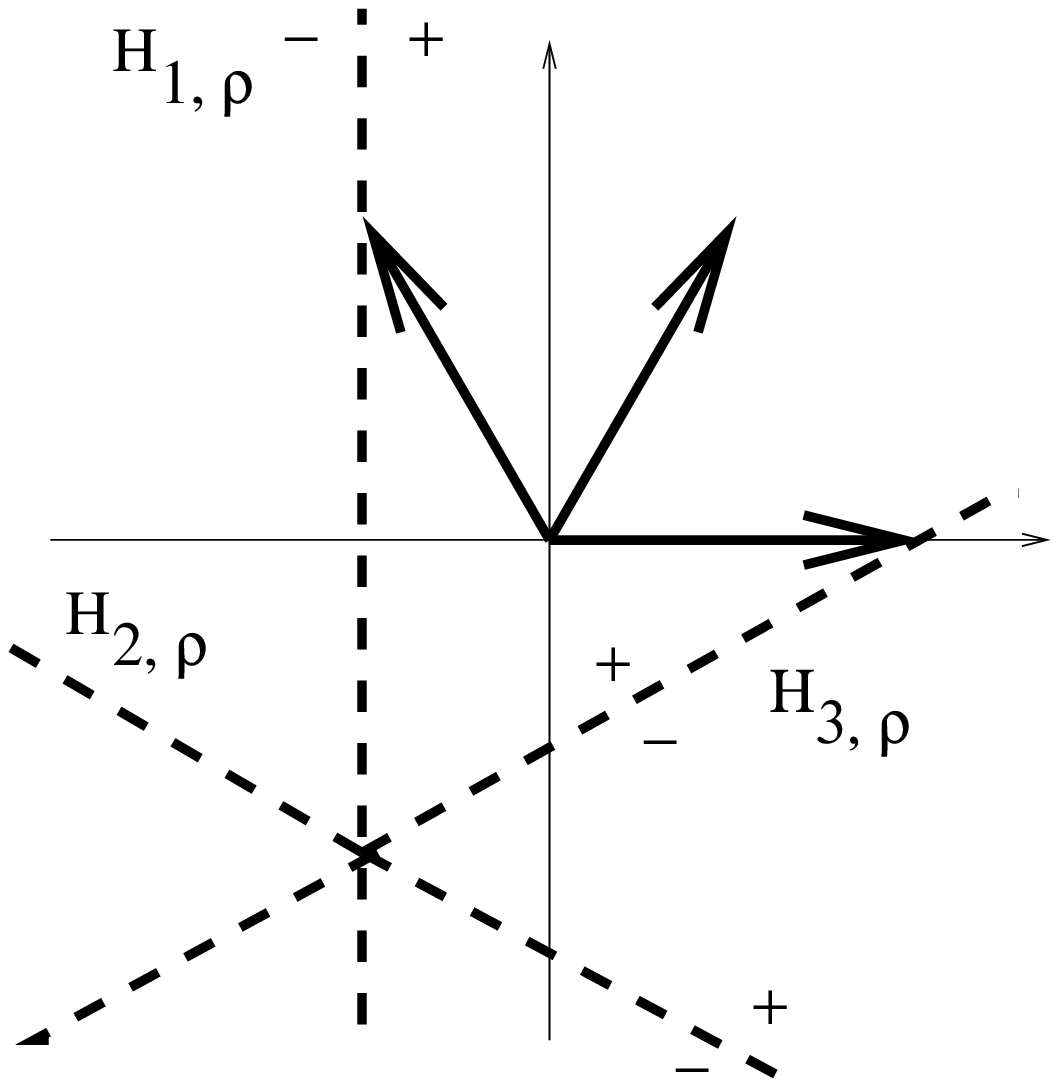,height=1.5in,width=1.5in} 
\label{cohomologychambers}
\epsfig{file=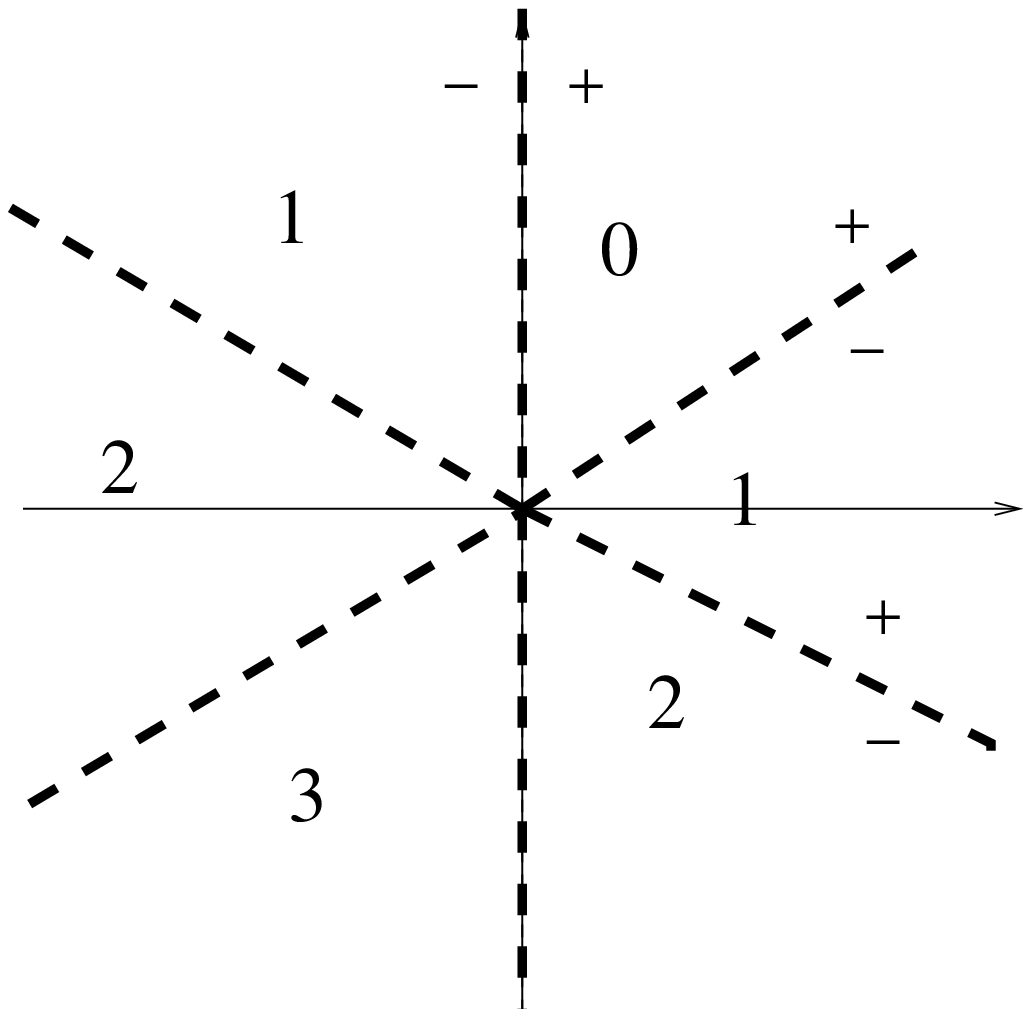,height=1.5in,width=1.5in} 
\label{asymptotic chambers}
\end{center}

For every chamber for the asymptotic cohomological functions, 
the number inscribed in it denotes the index of real
divisor classes in the chamber in question.
\end{example}

%% file: Estimates.tex
The content of this section is to establish 
the main technical tool in the proof of the continuity of asymptotic 
cohomological functions, the cohomological estimates based on the 
Mayer--Vietoris-type exact sequence of sheaves
\begin{eqnarray} \label{Mayer--Vietoris}
 && 0  \rightarrow  \sh{X}{D-\sum_{j=1}^{m}{A_j}} \rightarrow \sh{X}{D} \rightarrow \bigoplus_{1\leq i\leq m}\sh{A_i}{D} \rightarrow \\
&& \bigoplus_{1\leq i_1<i_2\leq m}\sh{A_{i_1}\cap A_{i_2}}{D} 
 \rightarrow \dots \rightarrow \bigoplus_{1\leq i_1< \dots < i_n \leq m}\sh{A_{i_1}\cap \dots \cap A_{i_n}}{D}\rightarrow 0\ , \nonumber
\end{eqnarray}  
where $D$ is an arbitrary integral divisor, $A, A_1,\dots , A_m$ are general  
very ample   divisors. The exactness of this  sequence 
 is established in Corollary \ref{long exact}. 
The cohomological  estimates obtained   via   this sequence will have 
a predominant role in the proof of the continuity of asymptotic cohomological
functions. 

We start out by establishing a local version of the sequence 
(\ref{Mayer--Vietoris}) under 
suitable general position hypotheses.
\begin{lem}\label{dual Koszul}
Let $R$ be a noetherian local ring, $n$ a nonnegative integer, 
$f_1,\dots, f_m\in R$ elements such that 
   \begin{enumerate}
     \item any $n$ element subset of $f_1,\dots , f_m$ forms a regular sequence in $R$,
     \item any $n+1$ elements from $f_1\dots , f_m$ generate $R$.
   \end{enumerate}
If $m<n$ then the complex 
\[ 
0\rightarrow (f_1\cdot\ldots\cdot f_m) \rightarrow R\rightarrow \bigoplus_{1\leq i\leq m} R/(f_i) \rightarrow
\dots \rightarrow R/(f_1,\dots, f_m) \rightarrow 0 \ 
\]
is exact.
If $m\geq n$ then 
\[
0\rightarrow (f_1\cdot\ldots\cdot f_m) \rightarrow R\rightarrow \bigoplus_{1\leq i\leq m} R/(f_i) \rightarrow 
\dots \rightarrow \bigoplus_{1\leq i_1<\dots <i_n\leq m} R/(f_{i_1},\dots , f_{i_n}) \rightarrow 0\ 
\]
is exact.
\end{lem}   

As this lemma is a dual version of Theorem 16.5 in \cite{Matsumura}, we only 
give an indication of the  proof. 

\begin{proof}
For any $1\leq i\leq m$ let $M(i)$ be the complex
\[
R\rightarrow R/(f_i)\rightarrow 0\ 
\]
with $R$ in degree $0$ and $R/(f_i)$ in degree $1$. 
Define $M^{(m)}$ as
$M^{(m)}\deq M(1)\otimes\dots\otimes M(m)$.
In the case $m\leq n$, $M^{(m)}$ is equal to 
\[
 R\rightarrow \bigoplus_{1\leq i\leq m} R/(f_i) \rightarrow
\dots \rightarrow R/(f_1,\dots, f_m) \rightarrow 0 \ ,
\]
while if $m> n$ then $M^{(m)}$ is 
\begin{eqnarray*}
&&  R\rightarrow \bigoplus_{1\leq i\leq m} R/(f_i) \rightarrow     \dots \rightarrow \bigoplus_{1\leq i_1<\dots <i_n\leq m} R/(f_{i_1},\dots , f_{i_n}) \rightarrow 0, 
\end{eqnarray*}
since all the quotients by ideals generated by at least $n+1$ of the $f_i$'s
are zero by assumption. By induction on $m$ and $\dim{R}$ one then proves 
that 
\[
H^{i}(M^{(m)}) = \begin{cases} (f_1\cdot\ldots\cdot f_m) & \mif i=0, \\
                               0 & \mif i\geq 1\ .
                \end{cases}
\]
\end{proof}

\begin{cor}\label{long exact}
Let $X$ be a pure $n$-dimensional  scheme of finite type over a field $k$, 
 $D$ an arbitrary Cartier divisor on $X$, $A_1,\dots, A_m$ effective  Cartier divisors on $X$, such that 
\begin{enumerate}
\item in the local rings of any point in $X$, the local equations of any collection of at most $n$ elements form
 a regular sequence  
\item the intersection of any $n+1$ of the $A_j$'s  is empty. 
\end{enumerate}
If  $m\leq n$ then the sequence
\[
 0\rightarrow \sh{X}{D-\sum_{j=1}^{m}{A_j}} \rightarrow \sh{X}{D} \rightarrow \bigoplus_{1\leq i\leq m}\sh{A_i}{D} \rightarrow \dots 
 \rightarrow \sh{A_1\cap\dots\cap A_m}{D} \rightarrow 0
\]
is exact, while if $m> n$ then
\begin{eqnarray*}
&& 0  \rightarrow  \sh{X}{D-\sum_{j=1}^{m}{A_j}} \rightarrow \sh{X}{D} \rightarrow \bigoplus_{1\leq i\leq m}\sh{A_i}{D} \rightarrow \\
&& \rightarrow \bigoplus_{1\leq i_1<i_2\leq m}\sh{A_{i_1}\cap A_{i_2}}{D}  \rightarrow \dots \rightarrow \bigoplus_{1\leq i_1< \dots < i_n \leq m}\sh{A_{i_1}\cap \dots A_{i_n}}{D}\rightarrow 0
\end{eqnarray*}                
is exact. Moreover, in $(1)$, it suffices to work with local rings at closed points.
\end{cor}

Our goal is to  provide  a cohomology estimate for certain differences of 
divisors which will serve as the basis of all further results.

\begin{lem} \label{cohomological estimate}
Let $C.  = C_{-1}\rightarrow C_0 \rightarrow C_1 \rightarrow \ldots $
be an exact sequence  of coherent sheaves on a proper scheme $X$ 
over a field. Then
\[
\left| \hh{i}{X}{C_{-1}}-\hh{i}{X}{C_0} \right| \leq 
2\cdot \sum_{k=1}^{i+1}{ \hh{i+1-k}{X}{C_k} }  .\]
\end{lem}

\begin{proof}
Break up the exact sequence $C.$ into a set of short exact sequences
\[  \ses{K_i}{C_i}{K_{i+1}} \ ,\]
where $0\leq i$ and $K_0=C_{-1}$. Observe that for every $r\geq 0$, 
a simple induction on $r$ shows that
\[
\hh{r}{X}{K_s}\ \leq\  \sum_{k=0}^{r}{ \hh{r-k}{X}{C_{s+k}}}\ .
\]
The  statement of the lemma then follows via 
\begin{eqnarray*}
\left| \hh{i}{X}{C_{-1}}-\hh{i}{X}{C_0} \right| & \leq & 
 \hh{i}{X}{K_1} +\hh{i-1}{X}{K_1} \\
& \leq & \sum_{k=1}^{i+1}{ \hh{i+1-k}{X}{C_k} } + 
\sum_{k=1}^{i}{ \hh{i-k}{X}{C_k}} \\
& \leq & 2\cdot \sum_{k=1}^{i+1}{ \hh{i+1-k}{X}{C_k} } \ .
\end{eqnarray*}

\end{proof}

\begin{cor}[Basic estimate] \label{basic estimate}
With notation as in Corollary $\ref{long exact}$, for any $i\geq 0$ and  $m\geq n$  
\begin{eqnarray*}
\left| \hh{i}{X}{D-\sum_{j=1}^{m}{A_j}} - \hh{i}{X}{D} \right|  & \leq &  
2\cdot  \sum_{k=1}^{i+1}{ \hh{i+1-k}{X}{\bigoplus_{1\leq i_1< \dots < i_k \leq m}\sh{A_{i_1}\cap \dots A_{i_k}}{D}   } } \ .
\end{eqnarray*}
\end{cor}

\begin{proof}
Apply Lemma \ref{cohomological estimate} to the long exact sequence in Corollary \ref{long exact}.
\end{proof}

%% file: Continuity.tex
The aim of this section is to prove the continuity of asymptotic 
cohomological functions. We will establish continuity in the following
form,  which  generalizes  the result  obtained for the 
volume function in \cite{pag},  Section 2.2.C. Our proof loosely 
follows the one given there. 
 
\begin{theorem}[Continuity of asymptotic cohomological functions]\label{main}
Let $X$ be an irreducible projective variety of dimension $n$. Then for all $0\leq i\leq n$,
\[ 
\widehat{h}^i: N^1(X)_{\QQ}\rightarrow {\RR}^{\geq 0}
\] 
defines a continuous  function on $N^1(X)_{\QQ}$ which is homogeneous of degree $n$, 
and satisfies the following Lipschitz-type estimate: there exists a constant $C$ 
such that for all pairs $\xi, \eta\in N^1(X)_{\QQ}$, one has 
\[
| \ha{i}{X}{\xi}-\ha{i}{X}{\eta} |\leq 
C\cdot \sum_{k=1}^{n}{ \zj{ \max\st{\norm{\xi},\norm{\eta}}}^{n-k}\cdot \norm{\xi-\eta}^k} 
\]
for some fixed norm $\|\ \|$.
\end{theorem}
\begin{remark}
As any two norms on a finite-dimensional real vector space are equivalent, 
it is indifferent which one we choose. However, the constant $C$ will 
certainly depend on it.
\end{remark}

\begin{corollary}\label{real version}
With notation as in the Theorem, the asymptotic cohomological 
functions $\widehat{h}^i$ extend uniquely to  continuous functions
\[
\widehat{h}^i: N^1(X)_{\RR}\rightarrow {\RR}^{\geq 0}\ ,
\]
which are homogeneous of degree $n$.
Moreover, they satisfy the same Lipschitz-type estimates as on 
rational classes in  the Theorem.
\end{corollary}

\begin{corollary} \label{cohomology and norm}
With notation as in the Theorem, there exists a positive 
constant $C$, such that for any real divisor class
 $\xi\in N^1(X)_{\RR}$  on $X$,  
\[ 
\ha{i}{X}{\xi} \leq C\cdot \norm{\xi}^n 
\]
for all $i$.
\end{corollary}

The proof of Theorem \ref{main} will come as the conclusion of a series 
of results,  
many  of which  are of somewhat  technical nature. Although some of them 
are probably well-known to experts, we will give proofs in cases when we 
were not able to find a suitable reference.  
Some statements which we need in the 
course of the proofs, but shed no further 
light to  our original problem will be relegated to Section 6.

The cohomological machinery developed in the previous section is only 
able to deal with pairs $\xi, \eta$  where
\[
\xi = \eta +\alpha,
\]
for some ample class  $\alpha\in N^1(X)_{\QQ}$, therefore we need 
some additional effort to pass to the general case.  This will rest upon 
 an  observation  about  normed vector  spaces (Proposition \ref{formal}).

We  establish the following 
statements,  which, when put together, will convey a  proof of 
Theorem \ref{main}.
As usual, $X$ will denote an irreducible complex projective variety of 
dimension $n$, all divisors are $\QQ$-Cartier unless explicitly mentioned 
otherwise.

\medskip\noindent
\textbf{Claim A} (Invariance with respect to numerical equivalence)
{\em
Let $D$ be an arbitrary, and $P$ a numerically trivial divisor. Then 
for all $i$, one has 
\[
\ha{i}{X}{D+P} = \ha{i}{X}{D}\ .
\]
}

\medskip\noindent
\textbf{Claim B} (Local uniform continuity in ample directions)
{\em 
Assume that the continuity of asymptotic cohomological functions 
holds for varieties of dimension less than $\dim (X)$. 
Let $D$ be an arbitrary, $A$ an  ample divisor on $X$. Fix an 
arbitrary  
norm on $N^1(X)_{\QQ}$. Then there exists a constant $C$ independent 
of $D$  for which
\[
\left| \ha{i}{X}{D-bA} - \ha{i}{X}{D} \right|\, \leq \,
C\cdot \sum_{k=1}^{n}{ \norm{D}^{n-k}\cdot b^k\norm{A}^k}\ .
\]
for every $b\geq1$.
}

\medskip\noindent
\textbf{Claim C} (Formal extension)
{\em
Let $V$ be an $r$-dimensional normed 
rational vector space, $A_1,\dots, A_r$ a basis 
for $V$, $f:V\rightarrow {\RR}^{\geq 0}$ a homogeneous function on $V$. 
Assume furthermore, that for every $1\leq i\leq r$ there exists a 
constant $C_i$ such that for all $D\in V$, and all natural numbers
$b\geq 1$,
\[
|f(D-bA_i)-f(D)| \,\leq\, C_i\cdot \sum_{k=1}^{n}{ \norm{D}^{n-k}\cdot b^k}\ .
\]
Then there exists a constant $C>0$, such that for every $D,D'\in V$ one has 
\[
|f(D)-f(D')| \,\leq\, C\cdot \sum_{k=1}^{n}{ 
\zj{\max\set{ \norm{D},\norm{D'}}}^{n-k}\cdot \norm{D-D'}^k}\ .
\]
}

These claims will appear below as  Proposition 
\ref{Claim B}, Proposition \ref{locunifamp}, and Proposition \ref{formal},
respectively. We now prove  Theorem \ref{main} based on these statements. 

\begin{proof}
According to Proposition \ref{Claim B} (Claim A), 
asymptotic cohomological functions are
invariant with respect to numerical equivalence of divisors, therefore they
are well-defined on $N^1(X)$. By their homogeneity property --- established
in Proposition
\ref{homogeneity} --- they are also well-defined on  $N^1(X)_{\QQ}$. Also,
asymptotic cohomological functions are homogeneous functions of degree 
$n=\dim X$ on $N^1(X)_{\QQ}$.

Observe that asymptotic cohomological functions are continuous on 
irreducible projective varieties of dimension $0$. By induction on the 
dimension of $X$,  Proposition \ref{locunifamp} (Claim B) then implies that, 
for any $\QQ$-divisor $D$ and any ample $\QQ$-divisor $A$, 
\[
\left| \ha{i}{X}{D-bA} - \ha{i}{X}{D} \right|\, \leq \,
C_A \cdot \sum_{k=1}^{n}{ \norm{D}^{n-k}\cdot b^k\cdot\norm{A}^k}
\]
for some fixed constant $C_A$, which is independent of $D$ (but possibly 
depends on $A$, and  certainly depends on the chosen norm) for 
all natural numbers $b\geq 1$.
Let $A_1,\dots A_r$ be a basis for $N^1(X)_{\QQ}$ consisting of ample
divisors. Then  the  functions $\widehat{h}^i$ on the rational vector space 
$N^1(X)_{\QQ}$ satisfy the assumptions in Proposition \ref{formal} (Claim C), 
and the theorem  follows.
\end{proof}

%**************************************************************************

\subsection{Rational continuity of asymptotic cohomological functions}

The first step along the way is to establish a  weak version of the rational 
continuity property of the volume.

\begin{proposition}\label{cohomology of intersections}
Let $X$ be a reduced projective scheme  of pure dimension $n$,  
$L$ an arbitrary line bundle, $A_1,\dots ,A_k$ very ample divisors on $X$, 
where $k\leq \dim X$. Then for every  pair of
ordered $k$-tuples $(E_1,\dots, E_k)$, $(E_1',\dots ,E_k')$ with 
 $E_j,E_j'\in |A_j|$ $(1\leq j\leq k)$ in a dense open subset of 
$|A_1|\times\ldots\times |A_k|$, 
one has
\[
\hh{i}{E_1\cap\dots \cap E_k}{L|_{E_1\cap\dots \cap E_k}}= 
\hh{i}{E_1'\cap\dots \cap E_k'}{L|_{E_1'\cap\dots \cap E_k'}}\ ,
\]
and the intersections $E_1\cap\dots \cap E_k$, $E_1'\cap\dots \cap E_k'$
are reduced.
\end{proposition}

\begin{proof}
For every $1\leq j\leq k$ let $V_j= {\PP}(\HH{0}{X}{A_i}^{\vee})$ 
be the projective space parameterizing  elements in $|A_j|$, 
let $T_j= \set{(x,D)\,|\, x\in D}\subseteq X\times V_i$ 
denote the total space of the flat family 

\medskip
\xymatrix{  & & & & T_j \ar[d]^-{{\phi}_j} \ar@{} [r] |-{ \subseteq}&  X\times V_j & \\
            & & & & V_j 
}
\medskip 
\noindent
with ${\phi}_j^{-1}(v)\cap \zj{X\times\st{v}}$ being the divisor in $X$ corresponding to $v\in V_j$. Consider the family 

\medskip
\xymatrix{  & & &  T_1{\times}_X\ldots{\times}_X T_k   \ar[d]^-{{\phi}_1\times\ldots\times {\phi}_k=\psi} \ar@{} [r] |-{ \subseteq} &  X\times (V_1\times\ldots\times V_k) & \\
            & & &  V_1\times\ldots\times V_k  &
}
\medskip
\noindent
which parametrizes  ordered $k$-tuples  of divisors $(E_1,\dots , E_k)$ with $E_j\in|A_j|$ for $1\leq j\leq k$ together with a specified point $x\in\bigcap_{j=1}^{k}{E_j}$. 
As $k\leq \dim X$, the map $\psi={\phi}_1\times\ldots\times {\phi}_k$ is surjective. 
Since $V_1\times\ldots\times V_k$ is integral, by generic flatness there exists a dense open set
\[
U\subseteq V_1\times\ldots\times V_k
\]
such that 
\[
\psi|_{{\psi}^{-1}(U)} : {\psi}^{-1}(U)\rightarrow U
\]
is flat. The map $\psi$ over $U$ is a flat family whose fibres are closed subschemes of the form $E_1\cap\ldots\cap E_k$, $E_j\in |A_j|$ for $1\leq j\leq k$.

By possibly shrinking $U$ one can arrange via  a Bertini-type argument that  
all the subschemes $E_1\cap\ldots\cap E_k\subseteq X$ that are the fibres of
$\psi$ over $U$ are actually (geometrically) reduced (this amounts to showing
that the set of points in the base over which the fibres are geometrically 
reduced is constructible, dense, and has maximal dimension).

 Consider the line bundle
\[
{\mathcal L} = i^*p_1^*L\ ,
\]
where
\[ 
i:  T_1{\times}_X\ldots{\times}_X T_k  \rightarrow X\times (V_1\times\ldots\times V_k) 
\]
is the inclusion map and 
\[
p_1 : X\times (V_1\times\ldots\times V_k) \rightarrow X
\]
is the first projection. For the line bundle ${\mathcal L}$ one has
\[ 
{\mathcal L}|_{{\psi}^{-1}(u)}\simeq L|_{E_1\cap\ldots\cap E_k}\ ,
\]
where  ${\psi}^{-1}(u)\cap X\times\st{u} = E_1\cap\ldots \cap E_k$ inside $X$. Then the statement of the lemma follows from the semi continuity theorem \cite{Hartshorne}, III.12. applied to ${\mathcal L}$ over the integral base $U$.

\end{proof}

\begin{corollary}
With notation as above, one has 
\[
\ha{i}{E_1\cap\dots\cap E_k}{L|_{E_1\cap\dots\cap E_k}} =
\ha{i}{E_1'\cap\dots\cap E_k'}{L|_{E_1'\cap\dots\cap E_k'}}
\]
for all $i$, $1\leq k\leq n$ and for every  pair of
ordered $k$-tuples $(E_1,\dots, E_k)$, $(E_1',\dots ,E_k')$ with 
 $E_j,E_j'\in |A_j|$ $(1\leq j\leq k)$ in a dense open subset of 
$|A_1|\times\ldots\times |A_k|$, for all $1\leq j\leq n$.
\end{corollary}

In order to be able to use Corollary \ref{basic estimate} to prove 
our rational  continuity result, we need the following Bertini-type 
result on intersections of  divisors in a very ample linear system.

\begin{lem}\label{generalized Bertini}
Let $X$ be a reduced  projective variety of pure dimension $n$ over an 
algebraically closed field, $A$ a very ample divisor on $X$. 
Then for any $m>n$ and  any general $(E_1,\dots E_m)$ in $ |A|$,
 we have that 
\begin{enumerate}
\item in the local ring of any point of $X$, the local equations 
defining any  cardinality $n$ subset of $E_1,\dots ,E_m$ form a 
regular sequence
\item the intersection of any $n+1$ of the $E_i$'s is empty. 
\end{enumerate}
In a similar vein, for every $1\leq m\leq n$ and general ordered $m$-tuple
 $(E_1,\dots E_m)$ in $ |A|$, one has that in the local ring of any  
point of $X$, the local equations defining $E_1,\dots ,E_m$ form a 
regular sequence. Moreover, in each case, we have that the scheme-theoretic
intersections of the $E_j$'s are geometrically reduced and equidimensional. 
\end{lem}
\begin{proof}
The same principles of constructibility, fibration and flatness that 
were used  in the previous lemma imply that in order  to prove 
that the corresponding subset of the parameter space of ordered $m$-tuples 
of divisors is open and dense, it is enough to see that it is nonempty. 
We will   proceed by induction. In the case $m=1$ we will pick 
$E_1\in |A|$ which is geometrically reduced and equidimensional 
(we can actually pick $E_1$ to be irreducible).

Assume we have $E_1,\dots , E_m\in |A|$ as required. Take any reduced 
equidimensional  Cartier divisor $E_{m+1}$ that does not contain 
any of the associated primes  of the reduced equidimensional 
$(n-k)$-dimensional scheme-theoretic intersections 
\[
E_{i_1}\cap \dots\cap E_{i_k}
\]
where $1\leq i_1<\ldots< i_k\leq m$ and $1\leq k\leq n$. Then  on one hand, $E_{m+1}$ will avoid all  of the $n$-fold intersections. On the other hand, $E_{m+1}$ is not a zero divisor in any of  intersections $E_{i_1}\cap \dots\cap E_{i_{n-1}}$ (where $i_1,\dots ,\ i_{n-1}\leq m$) hence locally forms a regular sequence with them. This proves both statements.
\end{proof}

\begin{remark}
Although we imposed the hypothesis of equidimensionality on $X$ and 
only obtained an equidimensional condition in the conclusion, when $X$ is
irreducible and we consider generic intersections 
\[
E_{i_1}\cap\ldots\cap E_{i_k}
\]
with $k<n$ in Lemma \ref{generalized Bertini} (for $m>n$ or $m\leq n$), then
these intersections are even irreducible by \cite{Jou}, Corollary 6.7.
\end{remark}

\begin{lemma} \label{choice}
Let $X$ be an irreducible projective variety of dimension $n$, $A$ a very
ample, $D$ an arbitrary divisor on $X$. Then for every $r\geq n$, 
there exist divisors 
$E_1,\dots ,E_r\in |A|$, such that for every $1\leq k\leq n$ and every 
choice $1\leq j_1<\dots j_k\leq r$, the intersection
\[
E_{j_1}\cap\ldots\cap E_{j_k}
\] 
is reduced, irreducible if $k<n$ (zero dimensional, if $k=n$). Furthermore, 
for all $0\leq i\leq n$ and $m\geq 1$, the values
\[
\hh{i}{E_{j_1}\cap\ldots\cap E_{j_k}}{mD}
\]
are minimal in the family parametrizing $k$-fold intersections of divisors. 
\end{lemma}

\begin{proof}
 Pick the divisors 
$E_1,\dots , E_n\in |A|$ in such a way that 
for any $1\leq k\leq n$ and any $1\leq j_1<\dots < j_k\leq r$, the
intersection
\[
E_{j_1}\cap\dots \cap E_{j_k}
\]
is reduced, irreducible if $k<n$ ($0$-dimensional if $k=n$), 
and for all $m\geq 1$, the values of 
\[
\hh{i}{E_{j_1}\cap\dots\cap E_{j_k}}{mpD|_{E_{j_1}\cap\dots\cap E_{j_k}}}
\]
are all minimal, ie. they are equal to the  
value which is taken up over a nonempty open set
of the variety parametrizing $k$-fold intersections in Proposition
\ref{cohomology of intersections}. Note that Proposition \ref{cohomology
of intersections} is applied separately for each $m\geq 1$. Naturally, as 
$m$ grows, the open locus may shrink. However,  
the Baire category theorem implies that there is a nonempty intersection
of all the open loci. 
\end{proof}

\begin{prop}[Rational continuity]\label{Rational continuity}
\label{Claim A}
Let $X$ be an irreducible projective variety of dimension $n$ over $\CC$, 
$D$ an  arbitrary divisor, $A$ an  ample divisor, $i\geq 0$. Then
\[ 
\frac{1}{p^n} \left| \ha{i}{X}{pD-A} - \ha{i}{X}{pD} \right| \rightarrow 0
\]
as $p\rightarrow \infty$.
\end{prop}

\begin{proof}
We first  reduce to the case when $A$ is very ample.
Let $m_0A$ be a fixed large integral multiple of $A$, which is very ample. 
By  the homogeneity of the asymptotic cohomological functions we see that 
\[
\frac{1}{p^n} \left| \ha{i}{X}{pm_0D-m_0A} - \ha{i}{X}{pm_0D} \right| =
\frac{m_0^n}{p^n} \left| \ha{i}{X}{pD-A} - \ha{i}{X}{pD} \right| \ .
\]
The proposition for $A$ very ample then  implies that the right-hand side 
converges to $0$ as $p\rightarrow \infty$. Hence we can assume from 
now on that $A$ is very ample. 

Fix a positive integer $p$, and   divisors 
$E_1,\dots , E_n\in |A|$ as in Lemma \ref{choice}, such that 
 for all $m\geq 1$, the values of 
\[
\hh{i}{E_{j_1}\cap\dots\cap E_{j_k}}{mpD|_{E_{j_1}\cap\dots\cap E_{j_k}}}
\]
are all minimal. Next, fix a natural number $m\geq n$. 
Then by Lemma \ref{generalized Bertini} and Proposition 
\ref{cohomology of intersections}
 we can find general divisors $F_1,\dots , F_m\in |A|$ 
which  satisfy the conditions of Corollary \ref{basic estimate}. Furthermore,
by Lemma \ref{choice}, we can assume that for all $1\leq k\leq n$ and all 
$1\leq j_1<\dots <j_k \leq m$
\[
\hh{i}{F_{j_1}\cap\dots\cap F_{j_k}}{mpD} =
\hh{i}{E_1\cap\dots\cap E_k}{mpD}\ .
\]
By  putting $mpD$ in place of $D$ 
and $F_1,\dots , F_m$ in place of $A_1,\dots ,A_m$ in Corollary 
\ref{basic estimate}, we obtain that 
\begin{eqnarray*}
\left| \hh{i}{X}{mpD-mA} - \hh{i}{X}{mpD} \right|  & = &  
\left| \hh{i}{X}{mpD-\sum_{j=1}^{m}{F_j}} - \hh{i}{X}{mpD} \right|  \\
& \leq &  2\cdot \estf{\rtlsum{i+1-k}} \\
 & =  &  2\cdot\estf{ \binom{m}{k} \hh{i+1-k}{E_1\cap\ldots\cap E_k}{mpD}} \ .
\end{eqnarray*}
Let us divide  both sides by $\tfrac{m^n}{n!}$ and take upper limits. Note
that the choice of the $F$'s  depends on $m$, therefore it is crucial 
that they got replaced by the $E's$. 

Using  Corollary \ref{limsup and cohomology}, Lemma \ref{limsup},
and the homogeneity of asymptotic cohomological functions,  
  we arrive at 
\begin{eqnarray*}
\left| \ha{i}{X}{pD-A}  -  \hh{i}{X}{pD} \right|     
&\leq & C_n\cdot \estf{ \limsup_m\zj{
\frac{\hh{i+1-k}{E_1\cap\dots\cap E_k}{mpD}}{m^{n-k}/{(n-k)}!}}} \\ 
 & = & C_n\estf{  \ha{i+1-k}{E_1\cap\dots\cap E_k}{pD}} \\
& = & C_n\estf{  p^{n-k}\ha{i+1-k}{E_1\cap\dots\cap E_k}{D}} \ ,
\end{eqnarray*}
where $C_n$ is a positive constant only depending on $n$.
Upon dividing by $p^n$ we arrive at the statement of the proposition.
\end{proof}

%***************************************************************************

\subsection{Numerical invariance of asymptotic cohomological functions}

The  crucial ingredient in  the proof of the numerical invariance 
of asymptotic cohomological functions 
is the fact that 
numerically trivial divisors form a bounded family. The specific 
version of this theme which we employ  is formulated in 
Proposition \ref{boundedness}.
 
As a first step, we show that one can give uniform estimates on the 
cohomology of multiples
of divisors in families. Note that this statement does not follow from
the semicontinuity theorem  in \cite{Hartshorne}, Section III.12.

\begin{proposition}[Cohomology estimate in families]\label{cohomology in families}
Let $f:{\mathcal X}\rightarrow T$ be a projective map of  noetherian schemes, 
$\shf$ a coherent sheaf on $\shx$, $\shl$ an invertible sheaf  on $\shx$. 

Then there exists a positive constant 
$C$ depending only on $f,T,\shl,\shf$  for which 
\[
\hh{i}{{\mathcal X}_t}{{\shf}_t\otimes {\shl}_t^{\otimes m}}\,
 \leq \, C\cdot m^{\dim{{\mathcal X}_t}} 
\]
for all $m\geq 1$, $i\geq 0$ and all $t\in T$.
\end{proposition}

\begin{proof}
We will proceed by induction on the maximal fibre dimension of $f$ and 
noetherian induction on $T$. The case $\dim \shx=0$ is straightforward to check.

For the inductive step, we can assume 
that the base is reduced and irreducible, as the fibres are   
unaffected by nilpotents in the base,  and
we can deal with the irreducible components one at a time. 

Our  strategy is to find a non-empty open subset $U_1\subseteq T$ over which 
the proposition holds with a certain constant $C_1$. Starting from 
there, we can (by noetherian induction) construct a stratification of $T$ 
into finitely many irreducible subschemes  $U_q\subseteq T,1\leq q\leq r$, 
such that the proposition holds over 
$U_q$ with a constant $C_q$. Then we reach the desired conclusion 
by setting 
\[
C\deq \max \set{C_1,\dots , C_r}\ .
\]

Let $\eta$ denote the generic point of $T$ and consider ${\mathcal X}_{\eta}$,
 the fibre over the generic point. 
Since ${\mathcal X}_{\eta}$ is projective,  we are able to find 
  very ample Cartier divisors 
${\sha}_{\eta}, {\shb}_{\eta}$ on ${\shx}_{\eta}$ such that 
\[
{\shl}_{\eta} = {\OO}_{{\shx}_{\eta}}({\sha}_{\eta})\otimes  
{\OO}_{{\shx}_{\eta}}(-{\shb}_{\eta})\ ,
\]
and ${\sha}_{\eta},{\shb}_{\eta}$ have  the properties that 
\begin{enumerate}
\item none of ${\sha}_{\eta},{\shb}_{\eta}$ contain any of the 
associated subvarieties of ${\mathcal F}_{\eta}$ on ${\shx}_{\eta}$,
\item
the local equation of ${\shd}_{\eta}$ in ${\shx}_{\eta}$ is not a zero-divisor
in ${\sha}_{\eta},{\shb}_{\eta}\subseteq {\shx}_{\eta}$.
\end{enumerate}

This way,  we obtain the short exact sequences
\[
\ses{{\shf}_{\eta}\otimes {\shl}_{\eta}^{\otimes m} \otimes {\OO}_{{\shx}_{\eta}}(-{\shb}_{\eta})}
{{\shf}_{\eta}\otimes {\shl}_{\eta}^{\otimes (m+1)}}
{ {\shf}_{\eta}\otimes {\shl}_{\eta}^{(m+1)}|_{{\sha}_{\eta}}}
\]
\[
\ses{ 
{\shf}_{\eta}\otimes {\shl}_{\eta}^{\otimes m}\otimes {\OO}_{{\shx}_{\eta}}(-{\shb}_{\eta}) 
} 
{ 
{\shf}_{\eta}\otimes {\shl}_{\eta}^{\otimes m}
} 
{ 
{\shf}_{\eta}\otimes {\shl}_{\eta}^{\otimes m}|_{{\shb}_{\eta}}
}\ .
\]

By generic flatness and denominator chasing, it is possible to extend ${\sha}_{\eta}$ and 
${\shb}_{\eta}$ to $U$-ample $U$-flat divisors ${\sha}$, ${\shb}$ over a 
non-empty open neighbourhood $U\subseteq T$ of 
$\eta$ in such a way that ${\shf}|_U$ is $U$-flat,  the divisors 
${\sha}_u,{\shb}_u$ are very ample for every $u\in U$, 
\[
{\shl}_{u} = {\OO}_{{\shx}_u}({\sha}_{u})\otimes  {\OO}_{{\shx}_u}(-{\shb}_{u})
\]
and ${\sha}_u,{\shb}_u$  do not contain the associated primes of ${\shf}_u$. 
Here  $U$-flatness ensures that formation of the ideals of ${\sha},{\shb}$ respects
base change on $U$. Moreover, it 
follows from Proposition 9.4.2. in \cite{EGAIV} that (by possibly shrinking
the  open subset $U\subseteq T$)  the following sequences are exact for 
all $u\in U$:
\[
\ses{{\shf}_{u}\otimes {\shl}_{u}^{\otimes m}\otimes {\OO}_{{\shx}_u}(-{\shb}_{u})}
{{\shf}_{u}\otimes {\shl}_{u}^{\otimes (m+1)}}
{ {\shf}_{u}\otimes {\shl}_{u}^{\otimes (m+1)}|_{{\sha}_{u}}}
\]
\[
\ses{{\shf}_{u}\otimes {\shl}_{u}^{\otimes m}\otimes  {\OO}_{{\shx}_u}(-{\shb}_u)}
{{\shf}_u\otimes {\shl}_{u}^{\otimes m}}
{ {\shf}_{u}\otimes {\shl}_{u}^{\otimes m}|_{{\shb}_{u}}}\ .
\]
From the corresponding long exact sequences we obtain that
\begin{eqnarray*}
&& \left| \hh{i}{{\mathcal X}_u}{ {\shf}_u\otimes {\shl}_{u}^{\otimes (m+1)}} - 
 \hh{i}{{\mathcal X}_u}{ {\shf}_u\otimes {\shl}_{u}^{\otimes m}} \right| \leq  \\
&& \hh{i-1}{{\sha}_u}{{\shf}_u\otimes {\shl}_{u}^{\otimes (m+1)}|_{{\sha}_{u}}} +
  \hh{i-1}{{\shb}_u}{{\shf}_u\otimes  {\shl}_{u}^{\otimes (m+1)}|_{{\shb}_{u}}} + \\
&& \hh{i}{{\sha}_u}{{\shf}_u\otimes {\shl}_{u}^{\otimes (m+1)}|_{{\sha}_{u}}} +
  \hh{i}{{\shb}_u}{{\shf}_u\otimes  {\shl}_{u}^{\otimes (m+1)}|_{{\shb}_{u}}} \ .
\end{eqnarray*}
The schemes  $\sha , \shb$ are  over $U$,
whose fibre dimension is strictly less than the maximal fibre dimension of 
${\shx}_U\rightarrow U$, ${\shf}|_{\sha}$, 
${\shf}|_{\shb}$ are coherent sheaves on the respective schemes, and 
${\shl}|_{\sha}$, ${\shl}|_{\shb}$ are both relative ($U$-flat) Cartier
divisors on $\sha,\shb$, respectively. This latter fact follows from the $U$-flatness
of ${\sha},{\shb}$. Hence we can apply the induction
hypothesis to the projective families 
$ \sha \rightarrow U\ ,\ \shb\rightarrow U$, and  obtain that 
\[
\hh{j}{{\sha}_u}{{\shf}_u \otimes {\shl}^{\otimes (m+1)}|_{{\sha}_u}} 
\leq C'm^{\dim{{\mathcal X}_u}-1}
\]
for all $j\geq 0$, where the constant $C'$ is independent of $u$ (it only depends on
${\sha}_U\rightarrow U,{\shf}|_{{\sha}_U}$ and ${\shl}|_{{\sha}_U}$),  
and similarly 
\[
\hh{j}{{\shb}_u}{{\shf}_u \otimes {\shl}^{\otimes (m+1)}|_{{\sha}_u}} 
\leq C''m^{\dim{{\mathcal X}_u}-1}
\]
Consequently, 
\[
 \hh{i}{{\mathcal X}_u}{ {\shf}_u\otimes {\shl}_{{\shx}_{u}}^{\otimes m}} 
 \leq C\cdot m^{ \dim{ {\mathcal X}_u } }
\]
for all $i\geq 0$ and all $u\in U$, where the positive constant $C$ is again independent of $u$.
\end{proof}

The 
uniform behaviour of numerically trivial divisors enters the 
picture  in the form of a 
vanishing theorem of Fujita (see the reference in the proof).

\begin{proposition} \label{boundedness}
Let $X$ be an irreducible projective variety of dimension $n$. Then there exists a family

\medskip
\xymatrix{  & & & & V \ar[d]^-{\phi} \ar@{} [r] |-{ \subseteq}&  X\times T & \\
            & & & & T 
}
\medskip

\noindent
with $T$ a quasi-projective variety (not necessarily irreducible), $V\subseteq X\times T$  
a closed subscheme and  $\phi$ flat, together with   a 
very ample divisor $A$ on $X$, such that 
\begin{enumerate}
\item $A+N$ is very ample for every numerically trivial divisor $N$ on $X$, 
\item if $D\in |A+N|$ for some $N{\sim}_{num}\, 0$, then $D=V_t$ for some $t\in T$.
\end{enumerate}
\end{proposition}

\begin{proof}
We start by showing that
 there exists a very ample line bundle $A$ on $X$ such that $A+N$ 
is very ample for every numerically trivial divisor $N$.

According to Fujita's vanishing 
theorem (\cite{pag} Section 1.4.D.), for any fixed ample 
divisor $B$ there exists $m_0>0$ such that 
\[
\HH{i}{X}{{\OO}_X(mB+E)}=0
\]
for all $m\geq m_0$, $i\geq 1$ and all nef divisors $E$. In particular, vanishing holds for every numerically trivial divisor. Take any very ample divisor $B$, let $m_0$ be as in Fujita's vanishing theorem. Consider $A'=(m_0+n)B+N$ where $N$ is an arbitrary numerically trivial divisor. Then by Fujita's theorem,
\[
\HH{i}{X}{{\OO}_X(A'-iB)}=\HH{i}{X}{{\OO}_X(N+(m_0+n-i)B)}=0
\] 
for all $i\geq 0$ hence $A$ is $0$-regular with respect to $B$. By Theorem 1.8.3 in 
\cite{pag}, $A'$ is globally generated. But then   $A'+B=(m_0+n+1)B+N$ is very ample. 
Observe that the coefficient of $B$ is independent of $N$ hence the choice 
\[
A\deq (m_0+n+1)B
\]
will satisfy the requirements for $A$.

To prove the Proposition, 
pick $A$ as above  to begin with. 
According to 1.4.36 in \cite{pag}, Section 1.4.D, 
there exists a scheme $Q$  of finite type over $\CC$,  
and a line bundle ${\mathcal L}$ on $X\times Q$, with the property that 
for every line bundle ${\OO}_X(A+N)$ with $N$ numerically trivial, 
there exists $q\in Q$ 
for which 
\[
{\OO}_{X\times {q}}(A+N) = {\mathcal L}|_{X\times\set{q}} \ .
\]
Let $p$,$\pi$  denote that projection maps from $Q\times X$ to $Q$ and  $X$, 
respectively.  We can  arrange by possibly twisting 
${\mathcal L}$ further by ${\pi}^*{\OO}_X(mA)$, that $R^jp_*\,{\mathcal L}=0$
for $j>0$, and that the natural map 
$\rho: p^*p_*\, {\mathcal L} \rightarrow {\mathcal L} $
is surjective. 

By the theorem on cohomology and base change (\cite{Hartshorne}, 
Section III.12)  
\[
{\mathcal E} \deq p_*\,{\mathcal L}
\]
is a vector bundle on $Q$ whose formation commutes with base change over $Q$. 
As ${\mathcal E}(q) = \HH{0}{X}{{\mathcal L}|_{X\times {q}}}$,
its projectivization ${\PP}(\mathcal E)$  will  parametrize the divisors
 $D\in |N+A|$, with $N$ numerically trivial.

The  universal divisor over ${\PP}(\mathcal E)$ is constructed as follows. 
The kernel $\shm$ of the  natural map  $\rho$ (which is surjective 
in our case) is  a vector bundle whose formation respects base change on $Q$. 
By restricting to the fibre of $p$ over $q$, one has an exact sequence 
\[
\ses{(M_X)_q}{\HH{0}{\set{q}\times X}{L}{\otimes}_{\CC} {\OO}_X}{L}
\]
where $L={\mathcal L}|_{\set{q}\times X}$ and
$(M_X)_q= \HH{0}{X}{{\shm}|_{\set{q}\times X}}$. Hence we see that 
\[
{\PP}_{\rm sub}(M_X) = \set{ (s,x)\,|\, x\in X,\ s\in \HH{0}{X}{L},\ s(x)=0}\ ,
\]
ie. $V={\PP}_{\rm sub}({\mathcal M})$ and $T={\PP}(\mathcal E)$ will satisfy
the requirements of the proposition.
\end{proof}

\begin{corollary} \label{boundedness for intersections}
With notation as above, fix a positive integer $1\leq k\leq n$, and a 
very ample divisor $A$ as in Proposition \ref{boundedness}.
Then  there exists a projective family 
\[
f: {\mathcal X} \rightarrow T'
\]
which parametrizes (possibly with repetitions) $k$-fold intersections 
$E_1\cap \dots\cap E_k$ of  all divisors such that  ${\OO}_X(E_j-A)$ is numerically 
trivial for $1\leq j\leq k$.
\end{corollary}
\begin{proof}
Consider the family $V\rightarrow T$ constructed in the Proposition, and 
proceed as in the proof of  Proposition \ref{cohomology of intersections}
using the same principles of flatness, fibration and constructibility. 
\end{proof}

\begin{lem} \label{uniform estimate}
With notation as above, let $P, A, D$ be divisors on $X$, with $P$  
numerically trivial, $A$  a very ample divisor 
for which   $A+P'$ is also  very ample for any choice of a numerically
trivial divisor $P'$, $D$ arbitrary. 
Fix $1\leq k\leq n$, and take $E_1,\dots E_k\in 
|A+P|$ such that $E_1\cap\ldots\cap E_k$ is geometrically reduced and of 
pure dimension $n-k$. 
Then there exists a positive constant $C$ only depending on $D$ and $A$, 
such that 
\[
\hh{i}{E_1\cap\dots \cap E_k}{mD} \leq C\cdot m^{n-k}
\]
for all $m\geq 1$.
\end{lem}

\begin{proof}
As we saw in Corollary \ref{boundedness for intersections}, 
there exists a projective family 
\[
f : \shx \rightarrow T
\]
which parametrizes $k$-fold intersections of divisors $E_1\cap \dots
\cap E_k$, with ${\OO}_X(E_j-A)$ numerically trivial for $1\leq j
\leq k$. In particular, all the intersections which occur in the lemma 
are fibres of this family. By restricting $f$ to the inverse image of those 
points in the base whose fibres have pure dimension $n-k$, we obtain an 
equidimensional family. Note that by \cite{EGAIV} 9.9.3. applied to 
\cite{EGAIV} 9.9.2 (1),(3) and (4), the set of points in the base 
whose fibres are of pure dimension $n-k$ is constructible, hence 
can be stratified by locally closed subsets. 

The claim of the lemma then follows from 
Proposition \ref{cohomology in families} 
applied to the line bundle 
\[
{\shl} \deq i^*pr_X^*{\OO}_X(D)\ ,
\]
where $pr_X$ is the projection map $\shx\rightarrow X$, and 
 $i$ is the inclusion of $\shx$ into $T\times X$.
\end{proof}

\begin{prop}[Invariance with respect to numerical equivalence]
\label{Claim B}
Let $X$ be an $n$-dimensional irreducible projective variety, 
$D$ an arbitrary,  $P$  a numerically trivial divisor on $X$. Then
\[  
\ha{i}{X}{D+P}=\ha{i}{X}{D}\ .
\]
In particular, the functions $\widehat{h}^i$ are well-defined on both 
$N^1(X)$, and $N^1(X)_{\QQ}$. 
\end{prop}

\begin{proof}
Fix a very ample divisor $A$ as in Lemma \ref{boundedness}, and fix a 
positive integer $p$. In particular,
for such  $A$, the divisor  $A+N$ is very ample for any 
numerically  trivial divisor $N$.
We want to estimate the difference
\[
\left| \hh{i}{X}{m(p(D+P)-A)} - \hh{i}{X}{mpD} \right|\ 
\]
in order to prove that 
\begin{equation} \label{numeric}
\frac{1}{p^n} \ha{i}{X}{p(D+P)-A}\rightarrow \ha{i}{X}{D}
\end{equation}
as $p\rightarrow \infty$. We know from Proposition \ref{Rational continuity} 
that 
\[
\frac{1}{p^n}\cdot \ha{i}{X}{p(D+P)-A} \rightarrow 
\ha{i}{X}{D+P}
\]
as $p\rightarrow \infty$. By the uniqueness of limits, (\ref{numeric})
then implies $\ha{i}{X}{D+P} = \ha{i}{X}{D}$.

Observe that 
\begin{eqnarray*}
\left| \hh{i}{X}{m(p(D+P)-A)} - \hh{i}{X}{mpD} \right|  
& = & \left| \hh{i}{X}{mpD-\sum_{j=1}^{m}{E_j^{(p)}}} - \hh{i}{X}{mpD} \right|
\end{eqnarray*} 
where $E_j^{(p)}\in |A-pP|$ for every $1\leq j\leq m$ are  general divisors. 
From Corollary \ref{basic estimate} we obtain
\begin{eqnarray*}
&& \left| \hh{i}{X}{mpD-\sum_{j=1}^{m}{E_j^{(p)}}} - 
 \hh{i}{X}{mpD} \right| \\  
&  \leq &  2\cdot \estf{\sum_{1\leq j_1<\ldots <j_k\leq m}{ \hh{i+1-k}{E_{j_1}^{(p)}\cap\ldots\cap E_{j_k}^{(p)}}{mpD}}} \ .
\end{eqnarray*}
Lemma \ref{uniform estimate} 
ensures the existence of a positive constant $C_{A,D}$, for which
\[
\hh{s}{E_{j_1}^{(p)}\cap \dots\cap E_{j_k}^{(p)}}{mpD} \leq C_{A,D}(mp)^{n-k}\ ,
\] 
and which  might depend on $D$ and  $A$, but is independent of $pP$, 
and the particular choices of the divisors  $E^{(p)}_{j_l}$. 
Therefore, 
\begin{eqnarray*}
\left| \hh{i}{X}{mpD-\sum_{j=1}^{m}{E_j^{(p)}}} - \hh{i}{X}{mpD} \right| &
\leq & C_{A,D}\estf{\sum_{1\leq j_1<\ldots <j_k\leq m}{ (mp)^{n-k}}} \\
& = &  C_{A,D}\estf{  \binom{m}{k}\cdot  (mp)^{n-k}}
\end{eqnarray*}
hence 
\[
\left| \hh{i}{X}{m(p(D+P)-A)} - \hh{i}{X}{mpD} \right|\ \leq C'm^np^{n-1}\ .
\]
After dividing by $\tfrac{m^n}{n!}$ and taking upper limits we arrive at 
\[ 
\left| \ha{i}{X}{p(D+P)-A}- \ha{i}{X}{pD} \right| \leq C''p^{n-1}
\]
which, upon dividing by $p^n$ gives the desired conclusion.
\end{proof}

%***************************************************************************

\subsection{Local uniform continuity in ample directions}

The subsection serves to prove the estimate of Theorem \ref{main} 
in ample directions. This forms the basis for the proof of the 
general case. 

\begin{prop} \label{locunifamp}
Let $X$ be an irreducible projective variety of dimension $n$, $D$ an 
arbitrary, $A$ an ample $\QQ$-divisor. 
Assume that the continuity of asymptotic cohomological functions holds 
for varieties of dimension less than $n$.
Fix an arbitrary  norm on 
$N^1(X)_{\QQ}$. Then there exists a constant $C_A$ independent of $D$ 
and such that
\[
\left| \ha{i}{X}{D-bA} - \ha{i}{X}{D} \right| \leq 
C_A\cdot \sum_{k=1}^{n}{ \norm{D}^{n-k}\cdot b^k\cdot\norm{A}^k}\ .
\]
for all integers $b\geq 1$.
\end{prop}

\begin{proof}
Both sides of the inequality are homogeneous 
of degree $n$,
therefore we can assume that both $D$ and $A$ are integral and $A$ 
is very ample. Also, assume without loss of generality that $\norm{A}=1$.

Fix divisors $E_1,\dots, E_n$, as in Lemma \ref{choice}, 
and choose   norms  
on $N^1(Y)_{\QQ}$  where $Y$ runs through all the possible intersections 
of the $E_j$'s, $1\leq j\leq n$, such that for any choice of $Y$ and 
any divisor $E$ on $X$ 
\[
\norm{E|_Y}\leq \norm{E}\ .
\]
As restriction of divisors induces a map of vector spaces
$N^1(X)_{\QQ} \rightarrow N^1(Y)_{\QQ}$, it is 
possible  to choose norms as required.

Next, fix a natural number $m\geq n$, and pick divisors $F_1,\dots , F_m$
as in Lemmas \ref{generalized Bertini} and \ref{choice}, such that for all 
choices $1\leq k \leq n$
and $1\leq j_1<\dots <j_k\leq m$ 
\[
\hh{i}{F_{j_1}\cap\ldots\cap F_{j_k}}{mD}= 
\hh{i}{E_1\cap\ldots\cap E_k}{mD} \ .
\]
Then by   Corollary \ref{basic estimate}, one has 
\begin{eqnarray*}
\left| \hh{i}{X}{mD-mA} - \hh{i}{X}{mD} \right| & = &  
\left| \hh{i}{X}{mD-\sum_{j=1}^{mb}{F_j}} - \hh{i}{X}{mD} \right|  \\
& \leq & 2\cdot \estf{ \sum_{1\leq l_1< \ldots < l_k\leq
mb}\hh{i-k+1}{F_{l_1}\cap\dots\cap F_{l_k}}{mD}}  \\
& = &   2\cdot  \estf{ \binom{mb}{k}\hh{i-k+1}{E_1\cap\dots\cap E_k}{mD}} \ .
\end{eqnarray*}
Let us divide  both sides by $\tfrac{m^n}{n!}$ and take upper limits. Then
by Corollary \ref{limsup and cohomology} and Lemma \ref{limsup}, one has 
\begin{eqnarray*}
\left| \ha{i}{X}{D-A} - \ha{i}{X}{D} \right| & \leq & 2\cdot 
\estf{ b^k\cdot
\limsup{ \frac{\hh{i-k+1}{E_1\cap\dots\cap E_k}{mD}}{m^{n-k}/{(n-k)}!}}}  \\
& = &  C_n\estf{b^k \ha{i-k+1}{E_1\cap\dots\cap E_k}{D}} \ .
\end{eqnarray*}
The  hypothesis of the Proposition 
--- in the form of Corollary \ref{cohomology and norm} --- implies that 
\[ 
\ha{i-k+1}{E_1\cap\dots \cap E_k}{mD} \leq C_k\cdot 
\norm{ D|_{E_1\cap \dots \cap E_k}}^{n-k} \ .
\]
We remark that  the $k$-fold intersections
of the divisors $E_j$ are reduced of the expected dimension, and either
irreducible or zero-dimensional, hence the induction hypothesis indeed 
applies to them. Hence 
\begin{eqnarray*}
\left| \ha{i}{X}{D-bA} - \ha{i}{X}{D} \right| 
 & \leq  &  C\cdot \estf{b^k \cdot C_k\cdot 
{\norm{D|_{E_1\cap\dots\cap E_k}}^{n-k}}} 
= C_A \sum_{k=1}^{n}{b^k  \cdot 
{\norm{D}^{n-k}}} 
\end{eqnarray*}
as required.
\end{proof}

%% file: Technical.tex
First we collect some simple properties of upper limits that we need. 

\begin{lem}[Properties of $\limsup$] \label{limsup}
Let $a_n,b_n,c_n$ be sequences of nonnegative real numbers. Then
\begin{enumerate}
\item $\limsup{(a_n+b_n)}\leq \limsup{a_n}+\limsup{b_n}$.
\item $\limsup{(a_nb_n)}\leq (\limsup{a_n})(\limsup{b_n})$.
\item If $|a_n-b_n| \leq c_n$ then 
\[
|\limsup{a_n}-\limsup{b_n}| \leq \limsup{c_n}\ .
\]
\end{enumerate}
\end{lem}

\begin{proof}
The first two statements are well-known. As for the third one, we observe that 
$ |a_n-b_n|\leq c_n$ is equivalent to
$a_n-b_n\leq c_n$  and $ b_n-a_n\leq c_n$,
that is, 
$a_n\leq  b_n+c_n$ and $ b_n\leq a_n+c_n$.
But (1) implies that 
\[ 
\limsup{a_n}\leq \limsup{b_n}+\limsup{c_n} \textrm{ and } \limsup{b_n}\leq \limsup{a_n}+\limsup{c_n}
\]
hence 
\[
|\limsup{a_n}-\limsup{b_n}| \leq \limsup{c_n}\ .
\]
\end{proof}

\begin{corollary}\label{limsup and cohomology}
Let $X$ be an irreducible projective variety on dimension $n$, $D,D'$ 
arbitrary divisors on $X$, such that 
\[
\left| \hh{i}{X}{{\OO}_X(mD)} - \hh{i}{X}{{\OO}_X(mD')} \right| \,\leq\, a_m
\]
for a sequence $a_m$ of nonnegative real numbers. Then
\[
\left| \ha{i}{X}{{\OO}_X(D)} - \ha{i}{X}{{\OO}_X(D')} \right| \,\leq \,
\limsup_m{\frac{a_m}{m^n/n!}}\ .
\]
\end{corollary}

Although a short exact sequence of sheaves does not give rise to 
anything similar to a long exact sequence for the asymptotic 
cohomological functions, one typical usage of the cohomology 
long exact sequence can to some extent be simulated. Namely, 
there is a substitute for the principle that the vanishings 
of certain cohomology groups give extra relations between others.

\begin{lemma}\label{support}
Let $X$ be an irreducible variety, 
\[
\ses{\sha}{\shb}{\mathcal C}
\]
be a short exact sequence of coherent sheaves on $X$, $L$ a line bundle
on $X$. Then if 
\[
\dim\supp\sha \leq n-1\ ,
\]
then 
\[
\limsup_{m}{ \frac{ \hh{i}{X}{\shb\otimes L^{\otimes m}}}{m^n/n!}} = 
\limsup_{m}{ \frac{ \hh{i}{X}{{\mathcal C}\otimes L^{\otimes m}}}{m^n/n!}}\ ,
\]
while if 
\[
\dim\supp{\mathcal C} \leq n-1\ ,
\]
then
\[
\limsup_{m}{ \frac{ \hh{i}{X}{\sha\otimes L^{\otimes m}}}{m^n/n!}} = 
\limsup_{m}{ \frac{ \hh{i}{X}{\shb\otimes L^{\otimes m}}}{m^n/n!}}\ .
\]
\end{lemma}

\begin{proof}
We will treat the case when $\dim\supp \sha \leq n-1$, the other case
can be dealt with the exact same way. Consider the short exact sequence 
\[
\ses{\sha\otimes {L^{\otimes m}}}{\shb\otimes {L^{\otimes m}}}
{{\mathcal C}\otimes {L^{\otimes m}}}\ .
\]
From the corresponding long exact sequence, we obtain
\[ \left| \hh{i}{X}{\shb\otimes {L^{\otimes m}}} - 
\hh{i}{X}{{\mathcal C}\otimes  {L^{\otimes m}}} \right| \leq 
 \hh{i-1}{X}{\sha\otimes {L^{\otimes m}}} + \hh{i}{X}{\sha\otimes
{L^{\otimes m}}}\ .
\]
As $\dim\supp \sha \leq n-1$, we have that 
\[
 \hh{i-1}{X}{\sha\otimes {L^{\otimes m}}} \leq C\cdot m^{n-1}
\textrm{ and }
 \hh{i}{X}{\sha\otimes {L^{\otimes m}}} \leq C\cdot m^{n-1}
\]
for some positive constant $C$. Then Corollary \ref{limsup and cohomology}
implies 
\[
\left| \limsup_{m}{ \frac{ \hh{i}{X}{\shb\otimes L^{\otimes m}}}{m^n/n!}} - 
\limsup_{m}{ \frac{ \hh{i}{X}{{\mathcal C}\otimes L^{\otimes m}}}{m^n/n!}} 
\right| \leq
\limsup_{m}{ \frac{ 2\cdot C\cdot m^{n-1}}{m^n/n!}} = 0\ .
\]
\end{proof}

The following result forms a part of the proof of the main theorem of 
this paper.

\begin{proposition}\label{formal}
Let $V$ be an $r$-dimensional normed 
rational vector space, $A_1,\dots, A_r$ a basis 
for $V$, $f:V\rightarrow {\RR}^{\geq 0}$ a homogeneous function on $V$. 
Assume furthermore, that for every $1\leq i\leq r$ there exists a 
constant $C_i$ such that for all $D\in V$, and all natural numbers
$b\geq 1$,
\[
|f(D-bA_i)-f(D)| \,\leq\, C_i\cdot \sum_{k=1}^{n}{ \norm{D}^{n-k}\cdot b^k}\ .
\]
Then there exists a constant $C>0$, such that for every $D,D'\in V$ one has 
\[
|f(D)-f(D')| \,\leq\, C\cdot \sum_{k=1}^{n}{ 
\zj{\max\set{ \norm{D},\norm{D'}}}^{n-k}\cdot \norm{D-D'}^k}\ .
\]
\end{proposition}

\begin{proof}
Assume that the given norm is the maximum norm with respect to the basis
$A_1,\dots , A_r$ of $V$. Let 
\[
D-D' = \sum_{j=1}^{r}{ b_jA_j}\ ,
\]
where by the homogeneity of both sides of the inequality in the proposition, 
we can assume that for all $1\leq j\leq r$, the coordinates $b_j$ are 
integers (not necessarily nonnegative).

We will  show that
\[
\left| f(D') - f(D) \right| \leq 
 C \cdot \sum_{k=1}^{n}{ \norm{D'}^{n-k}\cdot \norm{D-D'}^k}
\]
with some constant $C$ independent of $D$ and $D'$, 
from which the proposition follows by 
$\norm{D'}\leq\max\set{\norm{D},\norm{D'}}$.
Write 
\[
f\zj{D-\sum_{i=1}^{r}{ b_jA_j}} - f(D) \,=\,  
\sum_{l=1}^{r}{ \zj{f\zj{D- \sum_{j=1}^{l}{b_jA_j}} - 
f\zj{ D - \sum_{j=1}^{l-1}{b_jA_j}}}}\ ,
\]
ie. as a telescoping 
sum of terms, where each summand is a difference of two terms 
by a nonnegative multiple of a basis vector $A_j$. 
By the triangle inequality and 
the assumption of the lemma, 
\begin{eqnarray*}
&& \left| f\zj{D-\sum_{j=1}^{r}{b_jA_j}}-f\zj{D} \right| 
 \leq  \sum_{l=1}^{r}{\left| f\zj{ \zj{D-\sum_{j=1}^{l-1}{b_jA_j}}-b_lA_l}  -  f\zj{ D-\sum_{j=1}^{l-1}{b_jA_j}} \right| } \\
&  \leq & \sum_{l=1}^{r}{ C  \cdot \sum_{k=1}^{n}{ 
\max\set{\norm{D-\sum_{j=1}^{l-1}{b_jA_j}}^{n-k},
\norm{D-\sum_{j=1}^{l}{b_jA_j}}^{n-k}} \cdot \norm{b_lA_l}^k}}\ . 
\end{eqnarray*}

Applying the triangle inequality again, we obtain that 
\begin{eqnarray*}
\norm{ D-\sum_{j=1}^{l-1}{b_jA_j}} &  \leq  & 
\norm{D-\sum_{j=1}^{r}{b_jA_j}}+ \norm{ \sum_{j=l}^{r}{b_jA_j}} \ , \\ 
\norm{ D-\sum_{j=1}^{l}{b_jA_j}} &  \leq  & 
\norm{D-\sum_{j=1}^{r}{b_jA_j}}+ \norm{ \sum_{j=l+1}^{r}{b_jA_j}} \ ,
\end{eqnarray*}
therefore, by  Newton's binomial theorem  and collecting terms, one has
\begin{eqnarray*}
\left| f\zj{D-A} - f\zj{D} \right|  \leq  C\cdot \sum_{l=1}^{r}{ \sum_{k=1}^{n}{ \sum_{s=0}^{n-k}{ \binom{n-k}{s}\norm{D-\sum_{j=1}^{r}{b_jA_j}}^{n-k-s}\norm{ \sum_{j=l}^{r}{b_jA_j}}^s \cdot \norm{b_lA_l}^k }}} \ .
\end{eqnarray*}

Observe, that as we chose the maximum norm relative to the basis 
$A_1,\dots , A_r$ on $N^1(X)_{\QQ}$, one has  
\[
\norm{ \sum_{j=l}^{r}{b_jA_j}}^s =\max\st{b_l^s,\dots ,b_r^s}
\textrm{\ \  and \ \ }
b_l^k\cdot \max_{l\leq j\leq r}\st{b_j^s} \leq \max_{l\leq j\leq r} \st{b_j^{k+s}} \leq \norm{\sum_{j=1}^{r}{b_jA_j}}^{k+s}
\] 
for every $1\leq l\leq r$ and $1\leq k\leq n$ and $1\leq s\leq n-k$. 
This implies 
\begin{eqnarray*}
 \left| f\zj{D-A} - f\zj{D} \right|  
& \leq &  C\cdot  \sum_{l=1}^{r}{ \sum_{k=1}^{n}{ \sum_{s=0}^{n-k}{ \binom{n-k}{s}\norm{D-\sum_{j=1}^{r}{b_jA_j}}^{n-k-s} \cdot 
 \norm{\sum_{j=1}^{r}{b_jA_j}}^{k+s}}}} \\
& \leq & Crn(n!) \sum_{p=1}^{n}{ \norm{D-\sum_{j=1}^{r}{b_jA_j}}^{n-p}\cdot \norm{\sum_{j=1}^{r}{b_jA_j}}^{p}} \\
& = &  Crn(n!) \sum_{k=1}^{n}{ \norm{D'}^{n-k}\cdot \norm{D-D'}^{k}}\ .
\end{eqnarray*}
\end{proof}